\documentclass{article}
\usepackage[english]{babel}
\usepackage{lineno}
\usepackage{graphicx,multicol}
\usepackage{epic,eepic,epsfig}
\usepackage{caption}

\usepackage{amsmath,amsfonts,amsthm,amssymb}
\usepackage{pstricks,pstricks-add,pst-math,pst-xkey}

\usepackage{algorithm, algorithmic}

\usepackage{graphicx,pstricks,pst-node,amssymb}

\usepackage{graphicx}
\usepackage{amssymb}
\usepackage{color}
\newcommand{\jbj}[1]{{#1}}
\newcommand{\jbjr}[1]{{#1}}
\newcommand{\jbjb}[1]{{#1}}
\newcommand{\jh}[1]{{#1}}
\usepackage[geometry]{ifsym}
\newtheorem{theorem}{\hspace{5mm}Theorem}[section]
\newtheorem{proposition}[theorem]{\hspace{5mm}Proposition}
\newtheorem{definition}[theorem]{\hspace{5mm}Definition}
\newtheorem{lemma}[theorem]{\hspace{5mm}Lemma}
\newtheorem{corollary}[theorem]{\hspace{5mm}Corollary}
\newtheorem{conjecture}[theorem]{\hspace{5mm}Conjecture}
\newtheorem{problem}[theorem]{\hspace{5mm}Problem}
\newtheorem{question}[theorem]{\hspace{5mm}Question}

\newtheorem{claim}{Claim}
\newcommand{\DONOTTEX}[1]{}

\newcommand{\dom}{\rightarrow}
\newcommand{\pf}{\hspace{5mm}{\bf Proof: }}

\begin{document}
\bibliographystyle{plain}
\title{\bf Good orientations of 2T-graphs\thanks{\jbjb{The research of each author was supported by the Danish research council under grant number DFF-7014-00037B}}}
\author{J. Bang-Jensen\thanks{Department of Mathematics and Computer Science, University of Southern Denmark, Odense, Denmark (email: jbj@imada.sdu.dk); research supported by DISCO project, PICS-CNRS and the Danish research council}
    \and S. Bessy\thanks{LIRMM,
    Universit\'e de Montpellier, France (email:
    stephane.bessy@lirmm.fr); financial support from DISCO project, PICS-CNRS.}\and J. Huang\thanks{\jbj{Department of Mathematics and Statistics, University of Victoria, Victoria, B.C., Canada 
(email: huangj@uvic.ca); research supported by NSERC}}\and M. Kriesell\thanks{\jbj{Department of Mathematics, Technische Universit\"at Ilmenau, Germany (email: matthias.kriesell@tu-ilmenau.de)}}}
\maketitle

\begin{abstract}
  \setlength{\parindent}{0em}
  \setlength{\parskip}{1.5ex}
Graphs which contain $k$ edge-disjoint spanning trees have been characterized by
Tutte. Out-branchings and in-branchings are natural analogues of spanning trees 
for digraphs. Edmonds has shown that it can be decided in polynomial time whether 
a digraph contains $k$ arc-disjoint out-branchings or $k$ arc-disjoint 
in-branchings. Somewhat surprisingly, Thomassen proved that deciding whether 
a digraph contains a pair of arc-disjoint out-branching and in-branching is 
an NP-complete problem. This problem has since been studied for various classes of 
digraphs, giving rise to NP-completeness results as well as polynomial time 
solutions. 

\hspace{5mm} In this paper we study graphs which admit acyclic orientations that 
contain a pair of arc-disjoint out-branching and in-branching (such an orientation 
is called {\bf good}) and we focus on edge-minimal such graphs.  
A {\bf 2T-graph} is a graph whose edge set can be decomposed into two edge-disjoint
spanning trees. Vertex-minimal 2T-graphs with at least two vertices which are 
known as {\bf generic circuits} play an important role in rigidity theory for 
graphs. We prove that every generic circuit has a good orientation.
 Using this result we prove that if $G$ is 2T-graph whose vertex set has a partition $V_1,V_2,\ldots{},V_k$ so that each $V_i$  induces a generic circuit $G_i$ of $G$ and  the set of edges between different $G_i$'s form a matching in $G$, then $G$ has a good orientation.  We also obtain a characterization for  the case when the set of edges between different $G_i$'s form a {\bf double tree}, that is, if we contract each $G_i$ to one vertex, and delete parallel edges we obtain a tree. All our proofs are constructive and imply polynomial algorithms for finding the desired good orderings and the pairs of arc-disjoint branchings which certify that the orderings are good.

We also identify 
a structure which can be used to certify a 2T-graph which does not have 
a good orientation. 


  \end{abstract}

\section{Introduction}

We consider graphs and digraphs which may contain parallel edges and arcs 
respectively but no loops, and generally follow the terminology in \cite{bang2009}.
Our point of departure is the following theorem of Tutte, which
characterizes graphs that contain $k$ edge-disjoint spanning trees.

\begin{theorem} \cite{tutteJLMS36} \label{thm:tutte} 
A graph $G=(V,E)$ has $k$ edge-disjoint spanning trees if and only if, for every 
partition ${\cal F}$ of $V$, $e_{{\cal F}}\geq k(|{\cal F}|-1)$ where $e_{{\cal F}}$
is the number of edges with end vertices in different sets of ${\cal F}$.
\end{theorem}

Using matroid techniques, one can obtain a polynomial algorithm which either 
finds a collection of $k$ edge-disjoint spanning trees of a given graph $G$ or
a partition $\cal F$ for which $e_{{\cal F}} < k(|{\cal F}|-1)$ that shows no such
collection exists in $G$ (see e.g. \cite{recski1989}).

Let $D = (V,A)$ be a digraph and $r$ be a vertex of $D$. An {\bf
  out-branching} (respectively, {\bf in-branching)} in $D$ is a
spanning subdigraph $B^+_r$ (respectively, $B^-_r$) of $D$ in which
each vertex $v \neq r$ has precisely one entering (respectively,
leaving) arc and $r$ has no entering (respectively, leaving) arc. The
vertex $r$ is called the {\bf root} of $B^+_r$ (respectively,
$B^-_r$). It follows from definition that the arc set of an
out-branching (respectively, in-branching) of $D$ induces a spanning
tree in the underlying graph of $D$. It is also easy to see that $D$
has an out-branching $B^+_r$ (respectively, an in-branching $B^-_r$)
if and only if there is a directed path from $r$ to $v$ (respectively,
from $v$ to $r$) for every vertex $v$ of $D$.  A well-known result due
to Edmonds \cite{edmonds1973} shows that it can be decided in
polynomial time whether a digraph contains $k$ arc-disjoint
out-branchings or $k$ arc-disjoint in-branchings.

Somewhat surprisingly, Thomassen proved that the problem of deciding
whether a digraph contains a pair of arc-disjoint out-branching and
in-branching is NP-complete (see \cite{bangJCT51}). This problem has
since been studied for various classes of digraphs, giving rise to
NP-completeness results as well as polynomial time solutions
\cite{bangJCT51,bangJCT102,bangDAM161a,bangJGT42,bangC24,bangTCS438}.
In particular, it is proved in \cite{bangJGT42} that the problem is
polynomial time solvable for acyclic digraphs, and in
\cite{bangJCT102} that every 2-arc-strong locally semicomplete digraph
contains a pair of arc-disjoint out-branching and in-branching.

It turns out that acyclic digraphs which contain a pair of arc-disjoint 
out-branching and in-branching admit a nice characterization. 
Suppose that $D = (V,A)$ is an acyclic digraph and that $B^+_s,B^-_t$ are 
a pair of arc-disjoint out-branching and in-branching rooted at $s,t$ respectively 
in $D$. Then $s$ must be the unique vertex of in-degree zero and $t$ the unique 
vertex of out-degree zero in $D$. 
Let $X \subseteq V \setminus \{s\}$ and let $X^-$ denote the set of vertices with 
at least one out-neighbour in $X$. 
Since each vertex of $x \in X$ has an in-coming arc in $B^+_s$ and 
each vertex $x' \in X^-$ has an out-going arc in $B^-_t$, we must have 
\begin{equation}
\label{inoutacyc}
\sum_{x\in X^-}(d^+(x)-1)\geq |X|.
\end{equation}
The following theorem shows that these necessary conditions are also sufficient
for the digraph $D$ to have such a pair $B^+_s,B^-_t$.

\begin{theorem}\cite{bangJGT42} \label{acyclicin-outbr}
Let $D=(V,A)$ be an acyclic digraph in which $s$ is the unique vertex of in-degree 
zero and $t$ is the unique vertex of out degree zero. Then $D$ contains a pair of 
arc-disjoint out-branching and in-branching rooted at $s$ and $t$ respectively
if and only if (\ref{inoutacyc}) holds for every $X \subseteq V \setminus \{s\}$. 
Furthermore, there exists a polynomial algorithm which either finds the desired 
pair of branchings or a set $X$ which violates (\ref{inoutacyc}).
\end{theorem}

Every graph has an acyclic orientation. A natural way of obtaining an acyclic
orientation of a graph $G$ is to orient the edges according to a vertex ordering 
$\prec$ of $G$, that is, each edge $uv$ of $G$ is oriented from $u$ to $v$ if 
and only if $u \prec v$. In fact, every acyclic orientation of $G$ can be obtained 
in this way. Given a vertex ordering $\prec$ of $G$, we use $D_{\prec}$ to 
denote the acyclic orientation of $G$ resulting from $\prec$, and call $\prec$  
{\bf good} if $D_{\prec}$ contains a pair of arc-disjoint out-branching and 
in-branching. 
We also call an orientation $D$ of $G$ {\bf good} if $D = D_{\prec}$ for 
some good ordering $\prec$ of $G$. 
Thus a graph has a good ordering if and only if it has a good orientation.
We call such graphs {\bf good} graphs.
\jbjr{By Theorem \ref{acyclicin-outbr}, one can check in polynomial time whether 
a given ordering $\prec$ of $G$ is good and return a pair of arc-disjoint 
branchings in $D_{\prec}$ if $\prec$ is good.} However, no polynomial time 
recognition algorithm is known for graphs that have good orderings. 

An obvious necessary condition for a graph $G$ to have a good ordering is that
$G$ contains a pair of edge-disjoint spanning trees. This condition alone implies 
the existence of a pair of arc-disjoint out-branching and in-branching in an
orientation of $G$. But such an orientation may never be made acyclic for certain 
graphs, which \jh{means that $G$ does not have a good ordering}. On the other hand, to certify that 
a graph has a good ordering, it suffices to exhibit an acyclic orientation of $G$, 
often in the form of $D_{\prec}$, and show it contains a pair of arc-disjoint 
out-branching and in-branching. 
 
In this paper we focus on the study of edge-minimal graphs which have good 
orderings (or equivalently, good orientations). 

\begin{definition}
A graph $G = (V,E)$ is a {\bf 2T-graph} if $E$ is the union of two edge-disjoint 
spanning trees. 
\end{definition}

\jbj{Clearly, a graph has a good ordering if and only if it contains a spanning 
2T-graph which has a good ordering.} 
A 2T-graph on $n$ vertices has exactly $2n-2$ edges. 
The following theorem, due to Nash-Williams, implies a characterization of when 
a graph on $n$ vertices and $2n-2$ edges is a 2T-graph.
For a graph $G=(V,E)$ and $X\subseteq V$, the subgraph of $G$ {\bf induced} by 
$X$ is denoted by $G[X]$.

\begin{theorem}\cite{nashwilliamsJLMS39}
\label{thm:NWcover2T}
The edge set of a graph $G$ is the union of two forests if and only if 
\begin{equation}
  \label{sparse}
|E(G[X])|\leq 2|X|-2
\end{equation}
\noindent{}for every non-empty subset $X$ of $V$.
\end{theorem}

\jbj{
  \begin{corollary}
\label{NW2T}
A graph $G=(V,E)$ is a 2T-graph if and only if $|V| \geq 2$, $|E|=2|V|-2$,
 and (\ref{sparse}) holds.
    \end{corollary} 
  }

Generic circuits (see definition below) are important in rigidity theory for graphs.
A celebrated theorem of Laman \cite{lamanJEM4} implies that, for any graph $G$, 
the generic circuits are exactly the circuits of the so-called 
{\bf generic rigidity matroid} on the edges of $G$. 
Generic circuits have also been studies by Berg and Jord\'an \cite{bergJCT88},
who proved that every generic circuit is 2-connected and gave a full
characterization of 3-connected generic circuits (see Theorem \ref{thm:connelly}).

\begin{definition}
A graph $G=(V,E)$ is a {\bf generic circuit} if it satisfies the following
conditions:
\begin{enumerate}
\item[(i)] $|E|=2|V|-2>0$, and
\item [(ii)] $|E(G[X])|\leq 2|X|-3$, for every $X\subset V$ with 
                                   $2\leq |X| \leq |V|-1$.
\end{enumerate}
\end{definition}

Generic circuits are building blocks for 2T-graphs. According to Corollary 
\ref{NW2T}, each generic circuit is a 2T-graph on two or more vertices
with the property that no proper induced subgraph with two or more vertices is 
a 2T-graph. The only two-vertex generic circuit is the one having two parallel 
edges. Since no proper subgraph of a generic circuit is a generic circuit, every 
generic circuit on more than two vertices is a simple graph (i.e., containing no 
parallel edges). There is no generic circuit on three vertices and the only 
four-vertex generic circuit is $K_4$. \jbjb{The wheels\footnote{The {\bf wheel} $W_k$ is the graph that one obtains for a cycle of length $k$ by adding a new vertex and an edge from this to each vertex of the cycle.} $W_k$, $k\geq 4$ are all (3-connected) generic circuits.} Berg and Jord\'an \cite{bergJCT88} proved
that every 3-connected generic circuit can be reduced to $K_4$ by a series of
so called Henneberg moves (see definition below). We shall use this to prove
that every generic circuit has a good ordering.   

This paper is organized as follows. In Section \ref{2Tliftsec} we
begin with some preliminary results on generic circuits from
\cite{bergJCT88} and then prove a technical lemma that shows how to
lift a good orientation of a 2T-graph resulted from a Henneberg move
(Lemma \ref{lem:liftgoodor}). The lemma will be used in Section
\ref{OneGCsec} for the proof of a statement which implies \jbjb{that} every
generic circuit has a good ordering (Theorem \ref{thm:GChasgoodor}).
Section \ref{2TGCsec} is devoted to the study of the structure of
2T-graphs.  We show that every 2T-graph is built from generic circuits
and is reducible to a single vertex by a sequence of contractions of
generic circuits (Theorem \ref{thequotient}).  We also \jbjb{describe} a
polynomial algorithm which identifies all generic circuits of a
2T-graph (Theorem \ref{findallGC}). This implies that the problem of
deciding whether a graph is a disjoint union of generic circuits is
polynomial time solvable for 2T-graphs (Theorem
\ref{decomposeintoGC}). \jbjb{We also show that the problem is} NP-complete in general (Theorem
\ref{npc}). In Section \ref{good2T} we explore properties of 2T-graphs
which have good orderings and identify a forbidden structure for these
graphs (Theorem \ref{lem:pathconflict}). In Section
\ref{decomposeintoGCsec} we restrict our study on 2T-graphs which are
disjoint unions of generic circuits.  We prove that if the edges
connecting the different generic circuits form a matching, then one
can always produce a good ordering (Theorem \ref{matchingcase}) and we
also characterize when such an ordering exists if the graph reduces to
a double tree by contraction (Theorem \ref{thm:doubleTchar}).
Finally, in Section \ref{remarksec} we list some open problems and
show that the problem of finding a so called $(s,t)$-ordering of a
digraph is NP-complete (Theorem \ref{thm:storderNPC}).

\section{Lifting good orientations of a 2T-graph}\label{2Tliftsec}

\begin{definition}
  \label{Henneberg}
Let $G=(V,E)$ be a generic circuit, let $z$ be a vertex with three
distinct neighbours $u,v,w$.  A {\bf Henneberg move} from $z$ is the
operation that deletes $z$ and its three incident edges from $G$ and
adds precisely one of the the edges $uv,uw,vw$. A Henneberg move is
{\bf admissible} if the resulting graph, which we denote by
$G^{uv}_z$, where $uv$ is the edge we added to $G-z$, is a generic
circuit and a Henneberg move is {\bf feasible} if it is admissible and
$G^{uv}_z$ is a 3-connected graph.
\end{definition}

\begin{theorem} \cite{bergJCT88}
\label{thm:atleast3}
Let $G=(V,E)$ be a 3-connected generic circuit on $n\geq 5$
vertices. Then either $G$ has four distinct degree 3 vertices from
which we can perform an admissible Henneberg move, or $G$ has 3
pairwise non-adjacent vertices, each of degree 3, so that we can
perform a feasible Henneberg move from each of these.
\end{theorem}

\begin{theorem}\cite{bergJCT88}
\label{thm:connelly}
A 3-connected graph $G=(V,E)$ is a generic circuit if and only if $G$ can be reduced to (build from) $K_4$
by applying a series of feasible Henneberg moves (a series of Henneberg extensions\footnote{This is the inverse operation of a Henneberg move.}).
\end{theorem}

It is easy to see that if $z$ is a vertex of degree three in a 2T-graph $G$ then we
can obtain a new 2T-graph $G'$ by performing a Henneberg move from $z$.
The following lemma shows that when the three neighbours of the vertex $z$ that we 
remove in a Henneberg move are distinct, we can lift back a good orientation of 
$G'$ to a good orientation of $G$. 

\begin{figure}[h!t]
\begin{center}
\scalebox{0.9}{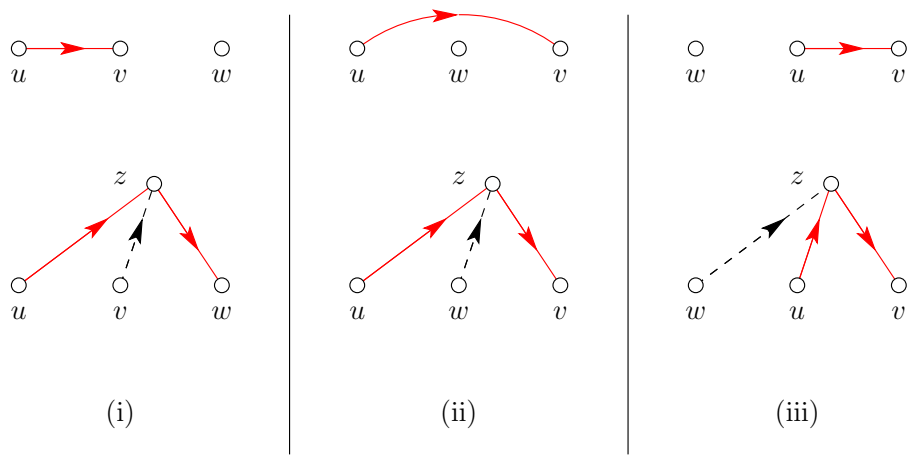}
\caption{How to lift a good ordering to a Henneberg extension as in Lemma \ref{lem:liftgoodor}. In-branchings are displayed solid, out-branchings are dashed.
The first line displays the three possible orders of the relevant vertices (increasing left to right) as they occur in the proof. The second line displays the ordering and the modification of the branchings in the extension.}
\label{F1}
\end{center}
\end{figure}

\begin{lemma}
\label{lem:liftgoodor}
Let $G$ be a 2T-graph on $n$ vertices and let $z$ be a vertex of degree 3 with three distinct neighbours $u,v,w$ from which we can perform an admissible Henneberg move to get $G^{uv}_z$. If  $G^{uv}_z$ is good, then also $G$ is good.
\end{lemma}
\pf Let ${\prec}'=(v_1,v_2,\ldots{},v_{n-1})$ be a good ordering of $G^{uv}_z$ and let $\tilde{B}^+_{v_1},\tilde{B}^-_{v_{n-1}}$ be arc-disjoint branchings of
$D_{\prec'}$.

Assume without loss of generality that $u=v_i$ and $v=v_j$ where $i<j$ (if this is 
not the case then consider the reverse ordering $\stackrel{\leftarrow}{\prec'}$ 
which is also good). Let $k\in [n-1]$ be the index of $w$ ($w=v_k$) and recall 
that $k\neq i,j$. Now there are 6 possible cases depending on the position of $w$ 
and which of the two branchings the arc $uv$ belongs to. 
In all cases we explain how to insert $z$ in the ordering ${\prec'}$ and update 
the two branchings which certifies that the new ordering $\prec$ is good.

\begin{itemize}
\item $uv$ is in $\tilde{B}^-_{v_{n-1}}$ and $j<k$. In this case we
  obtain ${\prec}$ from ${\prec'}$ by inserting $z$ anywhere between
  $v=v_j$ and $w=v_k$, replacing the arc $v_iv_j$ by the arcs
  $v_iz,zv_k$ and adding the arc $v_jz$ to $\tilde{B}^+_{v_1}$.
\item $uv$ is in $\tilde{B}^-_{v_{n-1}}$ and $i<k<j$. In this case we
  obtain ${\prec}$ from ${\prec'}$ by inserting $z$ anywhere between
  $w=v_k$ and $v=v_j$, replacing the arc $v_iv_j$ by the arcs
  $v_iz,zv_j$ and adding the arc $v_kz$ to $\tilde{B}^+_{v_1}$.
\item $uv$ is in $\tilde{B}^-_{v_{n-1}}$ and $k<i$. In this case we
  obtain ${\prec}$ from ${\prec'}$ by inserting $z$ anywhere between
  $u=v_i$ and $v=v_j$, replacing the arc $v_iv_j$ of
  $\tilde{B}^-_{v_{n-1}}$ by the arcs $v_iz,zv_j$ and adding the arc
  $v_kz$ to $\tilde{B}^+_{v_1}$.

\end{itemize}

The argument in the remaining three cases is obtained by considering   $\stackrel{\leftarrow}{{\prec'}}$ and noting that this switches the roles of the in- and out-branchings. \qed

\section{Generic circuits are all good}\label{OneGCsec}
%\section{Preliminaries on generic circuits}\label{prelimGC}

In this section we show that every generic circuit has a good ordering. In fact we prove a stronger statement on generic circuits \jbjb{which \jh{turns} out to be very} useful
in the study of 2T-graphs that have good orderings.

Let $H=(V,E)$ be 2-connected and let $\{u,v\}$ be a pair of non-adjacent vertices such that $H-\{u,v\}$ is not connected. Then there exists $X,Y\subset V$ such that $X\cap Y=\{u,v\}$, $X\cup Y=V$ and there are no edges between $X-Y$ and $Y-X$. A {\bf 2-separation} of $H$ along the cutset $\{u,v\}$ is the process which replaces $H$ by the two graphs $H[X]+e$ and $H[Y]+e$, where $e$ is a new edge connecting $u$ and $v$. It is easy to show the following.

\begin{lemma}\cite{bergJCT88} \label{lem:2sepgeneric}
Let $G=(V,E)$ is a generic circuit. Then $G$ is 2-connected. Moreover, if
$G-\{a,b\}$ is not connected, with connected components $X',Y'$, then $ab\not\in E$ and both of the graphs $G_1=G[X'\cup\{a,b\}]+ab$ and $G_2=G[Y'\cup\{a,b\}]+ab$ are generic circuits.
\end{lemma}

\iffalse
The inverse operation of 2-separation is that of a 2-sum: Given disjoint graphs $H_i=(V_i,E_i)$, $i=1,2$ and two prescribed edges $e_1=u_1v_1\in E_1$ and $e_2=u_2v_2\in E_2$, the {\bf 2-sum} $H_1\oplus{}H_2$ of $H_1,H_2$ along the pairs $u_1,u_2$ and $v_1,v_2$  is the graph we obtain from $H_1,H_2$ by deleting the edges $e_1,e_2$ and identifying $u_1$ with $u_2$ and $v_1$ with $v_2$.

\begin{lemma}\cite{bergJCT88}
Let $G_i=(V_i,E_i)$, $i=1,2$ be  generic circuits and let $u_iv_i\in E_i$, $i=1,2$ be edges. Then the 2-sum $G_1\oplus{}G_2$ along the pairs  $u_1,u_2$ and $v_1,v_2$ is a generic circuit.
\end{lemma}
\fi

\begin{theorem}
\label{thm:GChasgoodor}
Let $G=(V,E)$ be a generic circuit, let $s,t$ be distinct vertices of $G$ and 
let $e$ be an edge incident with at least one of  $s,t$.
Then \jbj{the following holds:
  \begin{enumerate}
  \item[(i)]  $G$ has a good ordering $\prec$ with corresponding branchings $B^+,B^-$ in which $s$ is the root of $B^+$, $t$ is the root of $B^-$  and $e$ belongs to $B^+$.
  \item[(ii)] $G$ has a good ordering $\prec$ with corresponding branchings $B^+,B^-$ in which $s$ is the root of $B^+$, $t$ is the root of $B^-$  and $e$ belongs to $B^-$.

  \end{enumerate}}
\end{theorem}

\begin{figure}[htbp]
\begin{center}
\scalebox{0.9}{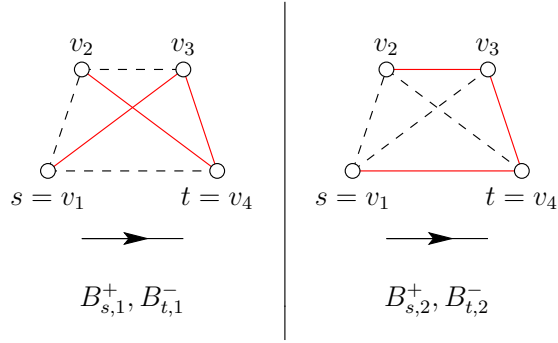}
\caption{Illustrating the base case of the inductive proof of Theorem
  \ref{thm:GChasgoodor}. All arcs are oriented from left to right, the
  prescribed edge is either $sv_3$ or $st$, and the picture shows that
  we may force $e$ to be in the (dashed) out-branching as well as in
  the (solid) in-branching.}
\label{F4}
\end{center}
\end{figure}

\pf The statement is clearly true when $G$ has two vertices. So assume
that $G$ has more than two vertices. The proof is by induction on $n$,
the number of vertices of $G$.  The smallest generic circuit on $n >
2$ vertices is $K_4$ and we prove that the statement holds for
$K_4$. By symmetry (reversing all arcs) it suffices to consider the
case when $e$ is incident with $s$ (see Figure~\ref{F4}).

It is possible to order the vertices of $K_4$ as $s=v_1,v_2,v_3,v_4=t$
such that $e \not =sv_2$, implying that $e=sv_3$ or $e=st$.  Let
$B^+_{s,1}$ and $B^+_{s,2}$ be the out-branchings at $s$ formed by the
arcs $sv_2,v_2v_3,st$ and of $sv_2,v_2t,sv_3$, respectively. The three
remaining edges form in-branchings $B^-_{t,1},B^-_{t,2}$ at $t$,
respectively. Since $st \in B^+_{s,1}$ and $sv_3 \in B^+_{s,2}$, we
find the desired branchings containing $e$ as in (i), (ii),
respectively.

\iffalse Consider an acyclic orientation of $K_4$ (it will be the
transitive tournament on 4 vertices) with the acyclic \jbj{ordering
  $v_1,v_2,v_3,v_4$ such that $s=v_1$, $t=v_4$ and the second end
  vertex of $e$ is either $v_3$ or $t$.}  \jbj{Suppose first that
  $e=st$. To get a pair of arc-disjoint branchings where $e$ belongs
  to $B^+$ we let $B^+$ consist of the arcs $sv_2,v_2v_3,st$ and let
  $B^-$ consist of the arcs $sv_3,v_2t,v_3t$. To get a pair of
  arc-disjoint branchings where $e$ belongs to $B^-$ we let $B^+$
  consist of the arcs $sv_2,sv_3,v_2t$ and $B^-$ consist of the arcs
  $st,v_2v_3,v_3t$.

  Suppose now that $e\neq st$.  If we want $e$ to be in $B^+$, then we
  let $B^+$ consist of the arcs $sv_2,sv_3,v_2t$ and let $B^-$ consist
  of the arcs $st,v_2v_3,v_3t$. If we want $e$ to be in $B^-$, then we
  let $B^+$ consist of the arcs $st,sv_2,v_2v_3$ and $B^-$ consist of
  the arcs $sv_3,v_2t,v_3t$.}\fi Assume below that $n>4$ and that the
statement holds for every generic circuit on at most $n-1$ vertices.

Suppose that $G$ is not 3-connected. Then it has a separating set
$\{u,v\}$ of size 2 (Recall that, by Lemma \ref{lem:2sepgeneric}, $G$
is 2-connected). Let $G_1,G_2$ be obtained from $G,u,v$ by
2-separation. Then each of $G_1,G_2$ are smaller generic circuits so
the theorem holds by induction for each of these. Note that $e$ is not
the edge $uv$ as this edge does not belong to $G$ by Lemma
\ref{lem:2sepgeneric}. \\

Suppose first that $s,t$ are both vertices of the same $G_i$, say w.l.o.g.
$G_1$. Then $e$ is also an edge of $G_1$ and there are two cases
depending on whether we want $e$ to belong to the out-branching or the
in-branching. We give the proof for the first case, the proof of the
later is analogous.

By induction there is a good ordering ${\prec}_1$ of $V(G_1)$ and
arc-disjoint branchings $B^+_{s,1},B^-_{t,1}$ so that $e$ belongs to
$B^+_{s,1}$. By interchanging the names of $u,v$ if necessary, we can
assume that the edge $uv$ is oriented from $u$ to $v$ in
$D({\prec}_1)$. Suppose first that the arc $u\rightarrow v$ is used in
$B^+_{s,1}$. By induction, by specifying the vertices $u,v$ as roots
and $e^+=uv$ as the edge, $G_2$ has a good ordering ${\prec}^+_2$ such
that $D({\prec}^+_2)$ has arc-disjoint branchings $B^+_{u,2},B^-_{v,2}$ where the
arc $uv$ is in $B^-_{v,2}$.  Now it is easy to check that
$B^+_s,B^-_t$ form a solution in $G$ if we let
\jbj{$A(B^+_s)=A(B^+_{s,1}-uv)\cup A(B^+_{u,2})$ and
  $A(B^-_t)=A(B^-_{t,1})\cup A(B^-_{v,2}-uv)$}. Here we used that
there is no edge between $u$ and $v$ in $G$, so $e$ is not the removed
arc above. The corresponding good ordering ${\prec}$ is obtained from
${\prec}_1,{\prec}^+_2$ by inserting all vertices of $V(G_2)-\{u,v\}$
just after $u$ in ${\prec}_1$. Suppose now that the arc $u\rightarrow
v$ is used in $B^-_{t,1}$. By induction, by specifying the vertices
$u,v$ as roots and $e^-=uv$ as the edge, $G_2$ has a good ordering
${\prec}^-_2$ and arc-disjoint branchings $B^+_{u,2},B^-_{v,2}$ such that the arc
$uv$ is in $B^+_{u,2}$. Again we obtain the solution in $G$ by
combining the two orderings and the branchings.  \jbj{By similar
  arguments we can show that there is also a good ordering such that
  the edge $e$ belongs to $B^-_t$.}

Suppose now that only one of the vertices $s,t$, say wlog. $s$ is a
vertex of $G_1$ and $t$ is in $G_2$. Note that this means that
$s,t\not\in \{u,v\}$.  Consider the graph $G_1$ with specification
$s,v,e$. By induction $G_1$ has a good orientation $D_1$ with
arc-disjoint branchings $B^+_{s,1},B^-_{v,1}$ so that $e$ belongs to
$B^+_{s,1}$. Note that, as $v$ is the root of the in-branching, the
edge $uv$ is oriented from $u$ to $v$ in $D_1$. If the arc $uv$
belongs to $B^+_{s,1}$, then we consider $G_2$ with specification
$u,t,uv$ where $uv$ should belong to the in-branching. By induction
there exists a good orientation $D_2$ of $G_2$ with arc-disjoint branchings
$B^+_{u,2},B^-_{t,2}$ such that the arc $uv$ is in \jbj{$B^-_{t,2}$}.
Now we obtain the desired acyclic orientation and arc-disjoint branchings by
setting $A(B^+_s)=A(B^+_{s,1}-uv)\cup A(B^+_{u,2})$ and
$A(B^-_t)=A(B^-_{v,1})\cup A(B^-_{t,2}-uv)$. To see that we do not
create any directed cycles by combining the acyclic orientations $D_1$
and $D_2$ it suffices to observe that $u$ has no arc entering in $D_2$
and $v$ has no arc leaving in $D_1$. If the arc $uv$ belongs to
$B^-_{v,1}$, then we consider $G_2$ with specification $u,t,uv$ where
$uv$ should belong to the out-branching. Again, by induction, there
exists an acyclic orientation $D_2$ of $G_2$ with good branchings and
combining the two orientations and the branchings as above we obtain
the desired acyclic orientation of $G$ and good in- and
out-branchings. \jbj{By similar arguments we can show that there is
  also a good ordering such that the edge $e$ belongs to $B^-_t$.}

It remains to consider the case when $G$ is 3-connected. By Theorem
\ref{thm:atleast3} there is an admissible Henneberg move $G\rightarrow
G^{uv}_z$ from a vertex $z\not\in\{s,t\}$ which is not incident with
$e$. Consider $G^{uv}_z$ with specification $s,t,e$, where $e$ should
belong to the out-branching. By induction there is an acyclic
orientation $D'$ of $G^{uv}_z$ and arc-disjoint branchings $B^+_s,B^-_t$ so that
$e$ is in $B^+_s$. Now apply Lemma \ref{lem:liftgoodor} to obtain an
acyclic orientation $D$ of $G$ in which $s$ is the root of an
out-branching $B^+$ which contains $e$ and $t$ is the root of an
in-branching $B^-$ which is arc-disjoint from $B^+$. The proof of the
case when $e$ must belong to $B^-_t$ is analogous. \qed

%then it follows from Theorem \ref{thm:connelly} that there is a series of $n-4$ feasible Henneberg moves leading from $G$ to $K_4$. However

\jbjr{

  We will see in Section \ref{decomposeintoGCsec} that Theorem
  \ref{thm:GChasgoodor} is very useful when studying good orderings of
  2T-graphs. The result below shows that it can also be applied to an
  infinite class of graphs which are not 2T-graphs.}

\begin{theorem}
  \label{4reg4con}
  Let $G$ be a $4$-regular $4$-connected graph in which every edge is
  on a triangle.  Then $G-\{e,f\}$ is a spanning generic circuit for
  any two disjoint edges $e,f$.  In particular, $G$ admits a good
  ordering.
\end{theorem}

{\bf Proof.}  Observe that $G$ is simple, as it is $4$-regular and
$4$-connected.  By Tutte's Theorem, $H:=G-\{e,f\}$ is a
2T-graph. Suppose, to the contrary, that it contains a 2T-graph $C$ as
a proper subgraph.  Then elementary counting shows that $C$ is an
induced subgraph of $G$ whose edge-neighborhood $N$ consists of
exactly four edges.  (In particular, neither $e$ nor $f$ connects two
vertices from $V(C)$.)  The endpoints of the edges from $N$ in $V(C)$
are pairwise distinct since $|V(C)| \geq 4$ and $G$ is $4$-connected.
Since $G-\{h,g\}$ is a 2T-graph for $h \not= g$ from $N$ we see that
$\overline{C}:=G-V(C)$ is a 2T-graph or consists of a single vertex
only.  If it is a 2T-graph then the endpoints of the edges of from $N$
in $V(\overline{C})$ are pairwise distinct, too, contradicting the
assumption that every edge is on at least one triangle. If, otherwise,
$\overline{C}$ consists of a single vertex only then it is incident
with both $e$ and $f$, contradicting the assumption that $e,f$ are
disjoint.
\hspace*{\fill}$\Box$

Thomassen conjectured that every $4$-connected line graph is
Hamiltonian \cite{thomassenJGT10}; more generally, Matthews and Sumner
conjectured that every $4$-connected claw-free graph (that is, a graph
without $K_{1,3}$ as an induced subgraph) is Hamiltonian
\cite{matthewsJGT8}. These conjectures are, indeed, equivalent
\cite{ryjacekJCT70}, and it suffices to consider $4$-connected line
graphs of cubic graphs \cite{kocholJCT78}.  Theorem \ref{4reg4con}
shows that such graphs have a spanning generic circuit (that is, a
spanning cycle in the rigidity matroid).

\jbj{
  \section{Structure of generic circuits in 2T-graphs}\label{2TGCsec}

Every 2T-graph $G$ on two or more vertices contains a generic circuit as 
an induced subgraph. Indeed, any minimal set $X$ with $|X| \geq 2$ and 
$|E(G[X])| = 2|X| -2$ induces a generic circuit in $G$. We say that $H$ is 
a {\bf generic circuit of} a graph $G$ if $H$ is a generic circuit and an 
induced subgraph of $G$.

\begin{proposition}
  \label{atmostoneincommon}
Let $G=(V,E)$ be a 2T-graph. Suppose that $G_1 = (V_1,E_1)$ and $G_2 =
(V_2,E_2)$ are distinct generic circuits of $G$. Then $|V_1 \cap
V_2|\leq 1$ and hence $|E_1 \cap E_2| = 0$. In the case when $|V_1
\cap V_2| = 0$, there are at most two edges between $G_1$ and $G_2$.
  \end{proposition}
  \pf Suppose to the contrary that $|V_1 \cap V_2| \geq 2$. Since
  $G_1$ and $G_2$ are generic circuits, $|E_1| = 2|V_1| - 2$ and
  $|E_2| = 2|V_2| - 2$.  Since $V_1 \cap V_2 \subset V_1$, we must
  have $|E_1 \cap E_2| = |E(G[V_1 \cap V_2])| \leq 2|V_1 \cap V_2| -
  3$. But then
  \begin{eqnarray*}
   |E(G[V_1\cup V_2])|&\geq & |E_1|+|E_2|-|E_1 \cap E_2|\\
                            &=& (2|V_1|-2)+(2|V_2|-2)-|E_1 \cap E_2|\\
                            &\geq& 2(|V_1|+|V_2|)-4-(2|V_1 \cap V_2| - 3)\\
                            &=& 2(|V_1| + |V_2| - |V_1 \cap V_2|)-1\\
                            &=& 2(|V_1 \cup V_2|) -1,\\
  \end{eqnarray*}
contradicting that $G$ is a 2T-graph and hence satisfies (\ref{sparse})
(see Corollary \ref{NW2T}). Hence $|V_1 \cap V_2| \leq 1$.

Suppose that $|V_1 \cap V_2| = 0$ (i.e., $G_1$ and $G_2$ have no vertex in common).
Let $k$ denote the number of edges between $G_1$ and $G_2$. Then
\begin{eqnarray*}
 \hspace{3cm}   k&=& |E(G[V_1\cup V_2])| - |E_1\cup E_2| \\
&\leq& (2|V_1 \cup V_2|-2) - (|E_1| + |E_2|)\\
&=& 2(|V_1|+|V_2|) -2 -(2|V_1|-2 + 2|V_2|-2)\\
&=& 2. \hspace{8cm}\qed \\
\end{eqnarray*}

\begin{proposition} \label{noedgebetween}
Let $r \geq 2$ and $G_i = (V_i,E_i)$ where $1 \leq i \leq r$ be generic circuits of
a 2T-graph $G = (V,E)$. Suppose that $|V_i \cap V_j| = 1$ if and only if 
$|i-j| = 1$. Then there is no edge with one end in $V_1 \setminus V_r$ 
and the other end in $V_r \setminus V_1$.
\end{proposition}
\pf Let $k$ be the number of edges each has one end in $V_1 \setminus V_r$ and the 
other end in $V_r \setminus V_1$. We prove $k = 0$ by induction on $r$.
When $r = 2$, 
\begin{eqnarray*}
 \hspace{1.5cm}   k&=& |E(G[V_1\cup V_2])| - |E_1\cup E_2| \\
&\leq& (2|V_1 \cup V_2|-2) - (|E_1| + |E_2|)\\
&=& 2(|V_1|+|V_2| -1) -2-(2|V_1|-2 + 2|V_2|-2)\\
&=& 0. \\
\end{eqnarray*}
Assume $r > 2$ and there is no edges with one end in $V_i \setminus V_j$ and 
the other end in $V_j \setminus V_i$ for all $1 \leq |i-j| \leq r-2$. By
assumption $|V_i \cap V_j| = 1$ if and only if $|i-j| = 1$ and in particular
$|V_i \cap V_j| = 0$ if $|i-j| > 1$. Hence
\begin{eqnarray*}
 k&=& |E(G[V_1\cup \cdots \cup V_r])| - |E_1\cup \cdots \cup E_r| \\
&\leq& (2|V_1 \cup \cdots \cup V_r|-2) - (|E_1| + \cdots + |E_r|)\\
&=& 2(|V_1|+ \cdots + |V_r| -(r-1)) -2-(2|V_1|-2 + \cdots + 2|V_r|-2)\\
&=& 0. \\
\end{eqnarray*} 
This completes the proof. 
\qed

  \begin{proposition}
    \label{hyperforest}
Let $G_i = (V_i,E_i)$ where $1 \leq i \leq r$ be the collection of generic 
circuits of a 2T-graph $G=(V,E)$ and let ${\cal G}=(V,\cal E)$ \jbjb{be the hypergraph} where 
${\cal E} = \{V_i:\ 1 \leq i \leq r\}$. Then $\cal G$ is a hyperforest.
  \end{proposition}
\pf Suppose to the contrary that $\cal G$ is not a hyperforest. Then there exist 
$V_{i_1},V_{i_2},\ldots{},V_{i_\ell}$ for some $\ell \geq 3$ such that
$|V_{i_j}\cap V_{i_k}|=1$ if and only if $|j-k|=1$ or $\ell-1$ and moreover,
the common vertices between the hyperedges on the hypercycle are pairwise distinct.
Thus 
$$\sum_{j=1}^{\ell}|V_{i_j}| = |V_{i_1}\cup{}V_{i_2}\cup\cdots{}\cup{}V_{i_\ell}| + \ell.$$  
By Proposition \ref{noedgebetween}, there is no edge with one end in
$V_{i_j}\setminus V_{i_k}$ and the other end in $V_{i_k} \setminus V_{i_j}$ for
all $j \neq k$. Hence 
    \begin{eqnarray*}
      |E(G[V_{i_1}\cup{}V_{i_2}\cup\cdots{}\cup{}V_{i_\ell}])|&=&\sum_{j=1}^{\ell}(2|V_{i_j}|-2)\\
      &=& 2|V_{i_1}\cup{}V_{i_2}\cup\cdots{}\cup{}V_{i_\ell}|,\\
    \end{eqnarray*}
  \noindent{}contradicting that $G$ is a 2T-graph and hence satisfies (\ref{sparse}). \qed 

Let $G = (V,E)$ be a 2T-graph and let ${\cal G}=(V,\cal E)$ be the hypergraph 
defined in Proposition \ref{hyperforest}. Two generic circuits of $G$ are
{\bf connected} if their vertex sets are in the same hypertree of ${\cal G}$. 
Not every vertex of $G$ needs to be in a generic circuit of $G$. 
A {\bf generic component} of $G$ is either a set consisting of a single vertex
which is not in \jbjb{any} generic circuit of $G$ or the union of a maximal set of 
connected generic circuits. A generic component is called {\bf trivial} if it
consists of a single vertex and {\bf non-trivial} otherwise. 
An edge of $G$ is {\bf external} if it is not contained
in any generic circuit. By Proposition \ref{noedgebetween} there is no
external edge in a generic component. Thus each generic component is 
a 2T-graph. Two generic components do not have a vertex in common.
A similar proof as for Proposition \ref{atmostoneincommon} shows that 
there can be at most two external edges between two generic components.
We summarize these properties below.

\begin{proposition}
\label{genericomponents}
Let $G = (V,E)$ be a 2T-graph. Then the following statements hold:
\begin{enumerate}
\item there is no external edge in a generic component;
\item each generic component is a 2T-graph;
\item two generic components are vertex-disjoint;
\item there are at most two external edges between two generic components.
\qed
\end{enumerate}
\end{proposition}  

Thus every 2T-graph $G$ partitions uniquely into pairwise vertex-disjoint
generic components. The {\bf quotient graph} $\tilde{G}$ of $G$ is the graph 
obtained from $G$ by contracting each generic component to a single vertex 
(and deleting loops resulted from the contractions). It follows from 
Proposition \ref{genericomponents} that every 2T-graph can be reduced to $K_1$ 
by successively taking quotients.

\begin{theorem}
\label{thequotient} 
Let $G$ be a 2T-graph. Then there is a sequence of 2T-graphs $G_0, G_1, \dots, G_k$
where $G_0 = G$, $G_k = K_1$, and $G_i = \tilde{G}_{i-1}$ for each 
$i = 1, 2, \dots, k$. In particular, $\tilde{G}$ is a 2T-graph.
\qed
\end{theorem}

\begin{theorem}
  \label{findallGC}
  There exists a polynomial algorithm $\cal A$ which given a 2T-graph $G=(V,E)$ as input finds the collection $G_1,G_2,\ldots{},G_r$, $r\geq 1$ of generic circuits of $G$.
\end{theorem}
\pf This follows from the fact that the subset system $M=(E,{\cal I})$ is a matroid, where $E'\subseteq E$ is in $\cal I$ precisely when $E'=\emptyset$ or $|E'|\leq 2|V(E')|-3$ holds, where $V(E')$ is the set of vertices spanned by the edges in $E'$. See \cite{bergLNCS2832} for a description of a polynomial  independence oracle. 
The circuits of $M$ are precisely the generic circuits of $G$. Recall from matroid theory that an element $e\in E$ belongs to a circuit of $M$ precisely when there exists a base of $M$ in $E-e$. Thus we can produce all the circuits by considering each edge $e\in E$ one at a time. If there is a base $B\subset E-e$, then $B\cup\{e\}$ contains a unique circuit $C_e$ which also contains $e$ and we can find $C_e$  in polynomial time by using independence tests in $M$. Since the generic circuits are edge-disjoint, by Proposition \ref{atmostoneincommon}, we will find all generic circuits by the process above.
\qed

\begin{corollary}
There exists a polynomial algorithm for deciding whether a 2T-graph $G$ is a generic circuit.
\end{corollary}

\begin{theorem}
  \label{decomposeintoGC}
  There exists a polynomial algorithm for deciding whether the vertex set of a 2T-graph $G=(V,E)$ decomposes into vertex disjoint generic circuits. Furthermore, if there is such a decomposition, then it is unique.
\end{theorem}

\pf We first use the algorithm $\cal A$ of Theorem \ref{findallGC} to find the set $G_1,G_2,\ldots{},G_r$ of generic circuits of $G$. 
If $r=1$ we are done as our decomposition consists of that generic circuit alone ($G$ is a generic circuit). So assume now that $r\geq 2$ and form the hypergraph $\cal G$ from $G_1,G_2,\ldots{},G_r$. Initialize
$H_1$ as the graph $G$ and ${\cal G}_1$ as the hypergraph $\cal G$.

By Proposition \ref{hyperforest}, ${\cal G}_1$ is a hyperforest and
hence, by Proposition \ref{atmostoneincommon} it has an edge which has
at most one vertex in common with the rest of the edges of ${\cal
  G}_1$. Let $G_{i_1}$ be a generic circuit corresponding to such an
edge. Note that, as $|V(G_{i_1})|\geq 2$ the generic circuit $G_{i_1}$
must be part of any decomposition of $V$ into generic circuits. Now
let $V_2=V-V(G_{i_1})$ and consider the induced subgraph $H_2=G[V_2]$
of $G$ and the hypergraph ${\cal G}_2=(V_2,{\cal E}_2)$ that we obtain
from ${\cal G}_1$ by deleting the vertices of $V(G_{i_1})$ as well as
every hyperedge that contains a vertex from $V(G_{i_1})$. If ${\cal
  G}_2$ has at least one edge, we can again find one which intersects
the rest of the edges in at most one vertex. Let $G_{i_2}$ denote the
corresponding generic circuit and add this to our collection. Form
$H_3,{\cal G}_3$ as above. Continuing this way we will either find the
desired decomposition or we reach a situation where the current
hypergraph ${\cal G}_k$ has at least one vertex but no edges. In this
case it follows from the fact that the generic circuits we have
removed so far are the only ones who could cover the corresponding
vertex sets that $G$ has no decomposition into generic circuits. As
the number, $r$, of generic circuits in $G$ is bounded by $|E|/2$
since generic circuits are edge-disjoint, the process above will
terminate in a polynomial number of steps and each step also take
polynomial time.\qed }
\medskip

The proof above made heavy use of the structure of generic circuits in
2T-graphs. For general graphs the situation is much worse.

\begin{theorem} \label{npc}
  It is NP-complete to decide if the vertex set of a graph admits a
  partition whose members induce generic circuits.
\end{theorem}

\jbjr{\pf Recall the problem {\sc exact
    cover by 3-sets} which is as follows: given a set $X$ with
  $|X|=3q$ for some integer $q$ and a collection ${\cal
    C}=Y_1,\ldots{},Y_k$ of 3-element subsets of $X$; does $\cal C$
  contain a collection of $q$ disjoint sets $Y_{i_1},\ldots{},Y_{i_q}$
  such that each element of $X$ is in exactly one of these sets?  {\sc
    exact cover by 3-sets} is NP-complete \cite[Page
    221]{garey1979}. Let {\sc exact cover by 4-sets} be the same
  problem as above, except that $|X|=4q$ and each set in $\cal C$ has
  size 4. It is easy to see that {\sc exact cover by 3-sets}
  polynomially reduces to {\sc exact cover by 4-sets}: Given an
  instance $X, \cal C$ of {\sc exact cover by 3-sets} we extend $X$ to
  $X'$ by adding $q$ new elements $z_1,z_2,\ldots{},z_q$ and construct
  ${\cal C}'$ by including the $q$ sets $Y\cup\{z_i\}$, $i\in [q]$ in
  ${\cal C}'$ for each set $Y\in \cal C$. It is easy to check that
  $X,{\cal C}$ is a yes-instance of {\sc exact cover by 3-sets} if and
  only if $X',{\cal C}'$ is a yes-instance of {\sc exact cover by
    4-sets} so the later problem is also NP-complete. Now given an
  instance $X',{\cal C}'$ of {\sc exact cover by 4-sets} we construct
  the graph $G$ as follows: the vertex set of $G$ consists of two sets
  $V_1$ and $V_2$. The set $V_1$ contains a vertex $v_x$ for each
  element $x\in X'$ and $V_2$ contains 4 vertices
  $u_{j,1},u_{j,2},u_{j,3},u_{j,4}$ for each set $Y_j\in {\cal C}'$ so
  that all these vertices are distinct. The edge set of $G$ is
  constructed as follows: for each $j\in [|{\cal C}'|]$ $E(G)$
  contains the edges of a $K_4$ on $u_{j,1},u_{j,2},u_{j,3},u_{j,4}$
  and for each set $Y_j\in {\cal C}'$ with $Y_j=\{x_1,x_2,x_3,x_4\}$
  $E(G)$ contains two copies of the edges
  $x_1u_{j,1},x_2u_{j,2},x_3u_{j,3}x_4u_{j,4}$. \jbjb{Clearly we can conctruct $G$ in polynomial time.}
  We claim that $G$ has
  a vertex partition into the vertex sets of disjoint generic circuits
  if and only if $X',{\cal C}'$ is a yes-instance of {\sc exact cover
    by 4-sets}.\\ Suppose first that $X',{\cal C}'$ is a yes-instance
  and let $Y_{i_1},\ldots{},Y_{i_q}$ be sets that form an exact cover
  of $X'$. For each $s\in [q]$ we include the 4 generic circuits of
  size 2 that connect the vertices $u_{s,1},u_{s,2},u_{s,3},u_{s,4}$
  to the vertices corresponding to $Y_{i_s}$ and for every other set
  $Y_j$ of ${\cal C}'$ (not in the exact cover) we include the generic
  circuit on the vertices $u_{j,1},u_{j,2},u_{j,3},u_{j,4}$. This
  gives a vertex partition of $V(G)$ into vertex sets of disjoint
  generic circuits. Suppose now that $G_1,\ldots{},G_p$ is a
  collection of vertex disjoint generic circuits such that $V(G)$ is
  the union of their vertex sets. Then we obtain the desired exact
  cover of $X'$ by including $Y_j\in {\cal C}'$ in the cover precisely
  when the generic circuit on $u_{j,1},u_{j,2},u_{j,3},u_{j,4}$ is not
  one of the $G_i$'s. Note that vertices of $V_1$ can only be covered
  by generic circuits of size 2 (parallel edges) so the sets we put in
  the cover will cover $X'$ and they will do so precisely once since
  $G_1,\ldots{},G_p$ covered each vertex of $G$ precisely once. Hence
  the chosen $Y_j$'s form an exact cover of $X'$.\qed

}

\section{Properties of good 2T-graphs} \label{good2T}

Let $G$ be a 2T-graph. For simplicity we shall call a generic circuit
of $G$ a {\bf circuit of} $G$. Recall from Section \ref{2TGCsec} that
each generic component of $G$ consists of either a single vertex or a
set of circuits that form a hypertree in ${\cal G}=(V,\cal E)$. We
call a generic component of $G$ a {\bf hyperpath} if its circuits
$G_1, G_2, \dots, G_k$ ($k \geq 1$) satisfy the property that for all
distinct $i, j$, $G_i$ and $G_j$ have a common vertex if and only if
$|i-j| = 1$. Note that the common vertices between circuits are
pairwise distinct and in particular, a generic component consisting of
one circuit is a hyperpath. We call $G$ {\bf linear} if every
non-trivial generic component of $G$ is a hyperpath.

\begin{proposition} \label{hyperpath}
Let $G$ be a 2T-graph which has one non-trivial generic component and
no trivial generic component. Then $G$ has a good ordering if and only
if $G$ is linear (i.e., $G$ is a hyperpath).
\end{proposition}
\pf Suppose that $G$ is a hyperpath formed by circuits $G_1, G_2,
\dots, G_k$.  For each $i = 1, 2, \dots, k-1$, let $v_i$ be the common
vertex of $G_i$ and $G_{i+1}$. Arbitrarily pick a vertex $v_0$ from
$G_1$ distinct from $v_1$ and a vertex $v_k$ from $G_k$ distinct from
$v_{k-1}$.  The assumption that $G$ is a hyperpath and the choice of
$v_0, v_k$ ensure that $v_0, v_1, \dots, v_k$ are pairwise distinct.
By Theorem \ref{thm:GChasgoodor}, each $G_i$ has a good ordering
$\prec_i$ that begins with $v_{i-1}$ and ends with $v_i$. It is easy
to see that the concatenation of these $k$ orderings gives a good
ordering of $G$.

On the other hand suppose that $G$ is not a hyperpath but has an
acyclic orientation with arc-disjoint branchings $B^+_s, B^-_t$. Since
$G$ is not a hyperpath, either there are three circuits intersecting
at the same vertex or there are three pairwise non-intersecting
circuits each intersecting with a fourth circuit. In either case, one
of the three circuits contains neither $s$ nor $t$. This would imply
that the arc sets of $B^+_s, B^-_t$ restricted to this circuit
contains a directed cycle, a contradiction to the fact that the
orientation of $G$ is acyclic.  \qed

\begin{proposition} \label{eachhyperpath}
If a 2T-graph $G$ has a good ordering, then $G$ is linear.
\end{proposition}
\pf Suppose that $D$ is a good orientation of $G$ with arc-disjoint
branchings $B^+_s, B^-_t$. Consider a non-trivial generic component
$H$ of $G$ and its orientation $D'$ induced by $D$ which is clearly
acyclic. Since $H$ has $2|V(H)|-2$ edges, $A(D') \cap A(B^+_s)$ and
$A(D') \cap A(B^-_t)$ induce arc-disjoint branchings in $D'$,
certifying that $D'$ is a good orientation of $H$. By Proposition \ref{hyperpath},
$H$ is is a hyperpath. Hence every non-trivial generic component of $G$ is 
a hyperpath and therefore $G$ is linear.  \qed
  
In view of Proposition \ref{eachhyperpath}, we only need to consider linear
2T-graphs for possible good orderings or good orientations. Suppose that $D$ is 
a good orientation of a 2T-graph $G$ with arc-disjoint branchings $B^+_s, B^-_t$. 
Let $H$ be a generic component of $G$. Then the proof of Proposition
\ref{eachhyperpath} shows that \jbjb{$D'=D[V(H)]$  is a good orientation of $H$  with 
arc-disjoint branchings $B^+_{s'}, B^-_{t'}$ which are the restrictions of $B^+_s, B^-_t$ to $V(H)$.}  We refer $s, t$ to as {\bf global} 
roots and $s',t'$ as to {\bf local} roots (of the corresponding branchings in $H$).
The {\bf external degree} of a vertex $x$ in $G$ is the number of external edges 
incident with $x$ in $G$ and the {\bf external degree} of $H$ is the sum
of external degrees of the vertices of $H$.

\begin{lemma} \label{localroots}
Let $G$ be a 2T-graph which has a good orientation with \jbjb{arc-disjoint} branchings $B^+, B^-$.
Then every non-trivial generic component has distinct local roots.
Suppose that $H, H'$ are generic components of $G$ and $xy$ is an external edge 
where $x, y$ are vertices in $H, H'$ respectively. 
\iffalse If $H$ does not contain a
global root and $x$ has external degree one (i.e., $xy$ is the only external edge 
incident with $x$), then one of the following holds:\fi
\jbjb{Then one of the following holds:}
\begin{itemize}
\item[(a)] $xy \in A(B^-)$ and $x$ is the local root of \jbjb{the in-branching  $B^-_H$ in $H$ which is the restriction of $B^-$ to $V(H)$};
\item[(b)] $xy \in A(B^+)$,  $y$ is the local root 
      of \jbjb{the out-branching  $B^+_{H'}$ in $H'$ which is the restriction of $B^+$ to $V(H')$} and if the external degree of $x$ is one, then $x$ is either the  root of $B^+$ (and hence the local
root of $B^+_H$  which is the restriction of $B^+$ to $V(H)$) or not a local root.
\item[(c)] $yx \in A(B^+)$ and $x$ is the local root of \jbjb{the out-branching  $B^+_H$ in $H$ which is the restriction of $B^+$ to $V(H)$}
\item[(d)] $yx \in A(B^-)$, \jh{$y$ is the local root of the in-branching $B^-_{H'}$ in $H'$ which is 
the restriction of $B^-$ and  if the external degree of $x$ is one, then $x$ is either the root of $B^-$ (and hence the local
root of $B^-_H$ which is the restriction of $B^-$ to $V(H)$) or not a local root.}
\end{itemize}
In particular, if the external degrees of $x, y$ are both one, and 
neither $H$ nor $H'$ contains a global root, then either $x$ is a local root
in $H$ or $y$ is a local root in $H'$ but not both. 
\end{lemma}

\jbjb{\pf Suppose that  $D$ is a good orientation of $G$ with arc-disjoint branchings $B^+, B^-$.
  Then, as we mentioned above, for every non-trivial generic component $H$ the restrictions of $B^+, B^-$ to $V(H)$
  form a pair of arc-disjoint branchings in $D[V(H)]$ and since $D$ is acyclic, the roots of these branchings must be distinct. Thus the first part of the lemma holds. This implies that  the digraph $\tilde{D}$ that we obtain by contracting each non-trivial generic component to one vertex is a good orientation of the quotient $\tilde{G}$ of $G$ and the digraphs $\tilde{B}^+, \tilde{B}^-$ that we obtain from $B^+, B^-$ via this contraction are arc-disjoint in- and out-branchings of $\tilde{D}$. As every vertex which is not the root of an in-branching (out-branching) has exactly one arc leaving it (entering it) this implies that if some arc $uv$ of $B^+$ ($B^-$) enters (leaves) a non-trivial
  generic component, then $v$ ($u$) is the local out-root (in-root) of that component. Now it is easy to see that (a)-(d) hold. The last claim is a direct consequence of these and the fact that $B^+$ and $B^-$ are arc-disjoint. \qed
}

We say that a subset $X\subset V$ with $2\leq |X|\leq |V|-2$ is {\bf
  pendant at $x$ in $G$} if all edges between $X$ and $V(G)-X$ are
incident with $x$. Note that $X$ is pendant \jbj{at $x$ in $G$ if an
  only if $V-X$ is pendant at $x$ }in $G$.

\begin{lemma}
  \label{lem:rootinpendant}
If $X$ is pendant at $x$ in a good 2T-graph $G$, then every good
orientation $D$ of $G$ will have $|X\cap\{s,t\}|=1$, where $s$ and $t$
are the roots of arc-disjoint branchings $B^+_s,B^-_t$ that certify
that $D$ is good. That is, $X$ contains precisely one global root.
\end{lemma}

\pf Let $B^+_s,B^-_t$ be a pair of arc-disjoint branchings that
certify that $D$ is good and suppose that none of $s,t$ are in $X$.
Let $z\in X-x$ (such a vertex exists as $|X|>1$).  As $X$ is pendant
in $x$ the $(s,z)$-path in $B^+_s$ passes through $x$ and the
$(z,t)$-path in $B^-_t$ also passes through $x$, but then $D$ contains
a directed cycle, contradicting that it is acyclic.  Since $V-X$ is
also pendant at $x$, we see that $|X\cap \{s,t\}|=1$ must hold.  \qed

Let $G$ be a 2T-graph.  Suppose that $H$ is a generic component of $G$ which is 
a hyperpath formed by circuits $G_1, G_2, \dots, G_k$. Then $H$ is called 
{\bf pendant} if one of following conditions holds:
\vspace{-3mm}
\begin{itemize} 
\item $V(H)$ is a pendant set in $G$;
\item all vertices of $G_1$ have external degree zero or all vertices of $G_k$ 
      have external degree zero.
\end{itemize}

\begin{corollary} \label{pendantcomponent}
If $H$ is a pendant generic component of a 2T-graph $G$, then $H$ must
contain a global root.
\end{corollary}
\pf If $H$ is the only generic \jbjb{component} in $G$, then clearly it
contains a global root. So assume that $G$ has at least two generic
components. We show that $V(H)$ contains a pendant set. If all
vertices of $G_1$ have external degree zero, then $H$ has at least two
circuits and $V(G_1)$ is a pendant set in $G$.  Similarly, if all
vertices of $G_k$ have external degree zero, then $V(G_k)$ is a
pendant set in $G$. In any case $V(H)$ contains a pendant set and
hence a global root by Lemma \ref{lem:rootinpendant}.  \qed

\begin{figure}[h!tbp]
\begin{center}
\scalebox{0.7}{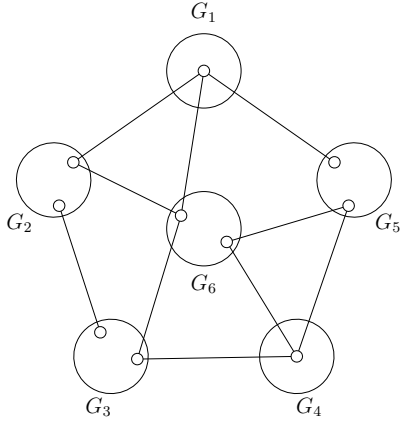}
\caption{\jbjr{Example of a 2T-graph $G$ whose vertex set is
    partitioned in circuits but which has no good ordering.  By
    Corollary \ref{pendantcomponent}, in any good orientation of $G$,
    the global roots $s,t$ are necessarily contained in $G_1$ and
    $G_4$. Now Lemma \ref{localroots} implies that the two
    vertices of attachment of $G_2,G_3$ must be local roots (of
    $G_2,G_3$, respectively) but not global roots in any good
    orientation. However, if two such local roots from distinct
    circuits are incident with only one external edge, then,
    by Lemma \ref{localroots}, these edges cannot be the same,
    implying that $G$ has no good ordering}}
\label{F5}
\end{center}
\end{figure}

\begin{corollary}
\label{cor:3pendant}
If $G$ contains three or more pairwise disjoint pendant subsets
$X_1,X_2,X_3$, then $G$ has no good orientation. In particular, a
2T-graph has a good ordering then it contains at most two pendant
generic components.
\end{corollary}
\pf  This follows immediately from Lemma~\ref{lem:rootinpendant} and Corollary 
\ref{pendantcomponent}.
\qed

 \begin{theorem}
   \label{lem:pathconflict}
   Suppose that there are vertex-disjoint circuits
   $G_{i_0},G_{i_1},\ldots{},G_{i_p}$, $p\geq 1$ of a 2T-graph $G$
   such that
   \begin{itemize}
   \item Each $G_{i_j}$ has external degree 3
   \item Some vertex $x_0\in V(G_{i_0})$ has external degree 2 and the
     third external edge goes between a vertex $y_0\in V(G_{i_0})-x_0$
     and a vertex $z_1\in V(G_{i_1})$
   \item Some vertex $x_p\in V(G_{i_p})$ has external degree 2 and the
     third external edge is adjacent to a vertex $y_p\in
     V(G_{i_p})-x_p$ and a vertex $z_{p-1}\in V(G_{i_{p-1}})$, where
     $z_{p-1}\neq x_0$ if $p=1$.
   \item For each $j\in [p-1]$ there is exactly one external edge
     between $V(G_{i_j})$ and $V(G_{i_{j+1}})$: $y_jz_{j+1}$ with
     $y_j\in V(G_j)$ and $z_{j+1}\in V(G_{j+1})$.
   \end{itemize}
   If $G$ has a good ordering $\prec$, then some vertex of
   $V(G_{i_0})\cup{}V(G_{i_1})\cup\ldots{}V(G_{i_p})$ is the first or
   the last vertex according to $\prec$ \jbjr{(that is, at least one
     of the global roots $s,t$ belongs to that vertex set).}
 \end{theorem}

\begin{figure}[h!tbp]
\begin{center}
\scalebox{0.7}{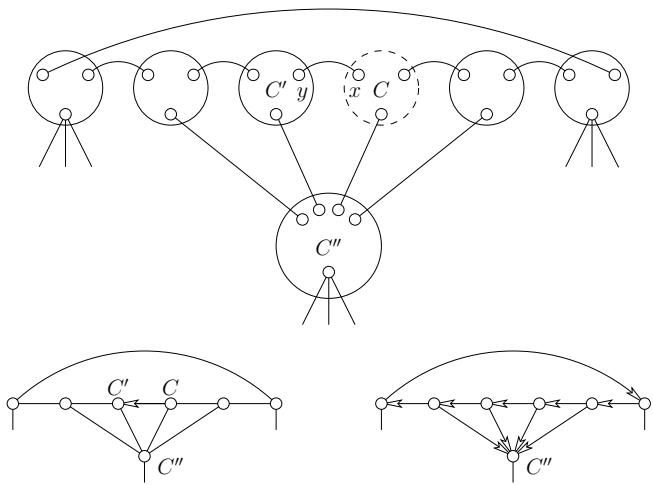}
\caption{\jbjr{The figure above shows part of a graph $G$ whose vertex
    set is partitioned in circuits together with all the external
    edges connecting them to other circuits.  Assume that we have a
    good ordering and that the seven circuits displayed in the
    configuration do not contain any of the two global roots.
    Consider the external edge $xy$ between $C$ and $C'$. By Lemma
    \ref{localroots}, exactly one of $x$ and $y$ must be a local
    root.  Say, w.l.o.g. that $x$ is a local root so $y$ cannot be a
    local root of $C'$. We encode this fact by a white arrow from $C$
    to $C'$ in the quotient graph (lower left figure). Now the other
    two vertices displayed in $C'$ must be its local roots, so that,
    following our drawing convention, we need to orient the remaining
    two edges incident with $C'$ in the quotient away from it.
    Processing this way all the six circuits in the upper row we get
    the lower right figure and deduce that, finally, there is no way
    to place two local roots in the circuit $C''$. Thus the conclusion
    is that if $G$ has a good ordering, then at least one of the
    global roots must be a vertex of one of the circuits in the upper
    part of the figure.  }}
\label{F6}
\end{center}
\end{figure}

\vspace{-5mm}

\pf Assume that $V(G_{i_0})\cup{}V(G_{i_1})\cup\ldots{}V(G_{i_p})$
does not contain any global root. The two local roots of $G_{i_0}$
are $x_0$ and $y_0$. So $z_1$ can not be a local root of $G_{i_1}$. 
Then $y_1$ is a local root of $G_{i_1}$ and $z_2$ is not a local root of
$G_{i_2}$. Following the argument, $z_p$ is not a local root of
$G_{i_p}$, but then it has only one local root, a contradiction.
\qed

 We call $G_{i_0},G_{i_1},\ldots{},G_{i_p}$ as above a {\bf conflict}
 of $G$.

 We say that two conflicts are {\bf disjoint} if no circuit is
 involved in both of them. %Clearly a type 1 conflict is disjoint from
 %every other conflict.

 The following is immediate from Corollary \ref{cor:3pendant} and
 Theorem \ref{lem:pathconflict}. For an example, see Figure \ref{F5}.

 \begin{corollary}
   \label{cor:atmost2conflicts}
   Let $G$ be a 2T-graph. If $G$ has 3 disjoint conflicts, then $G$
   has no good ordering.
 \end{corollary}

\begin{figure}[htbp]
\begin{center}
\scalebox{0.5}{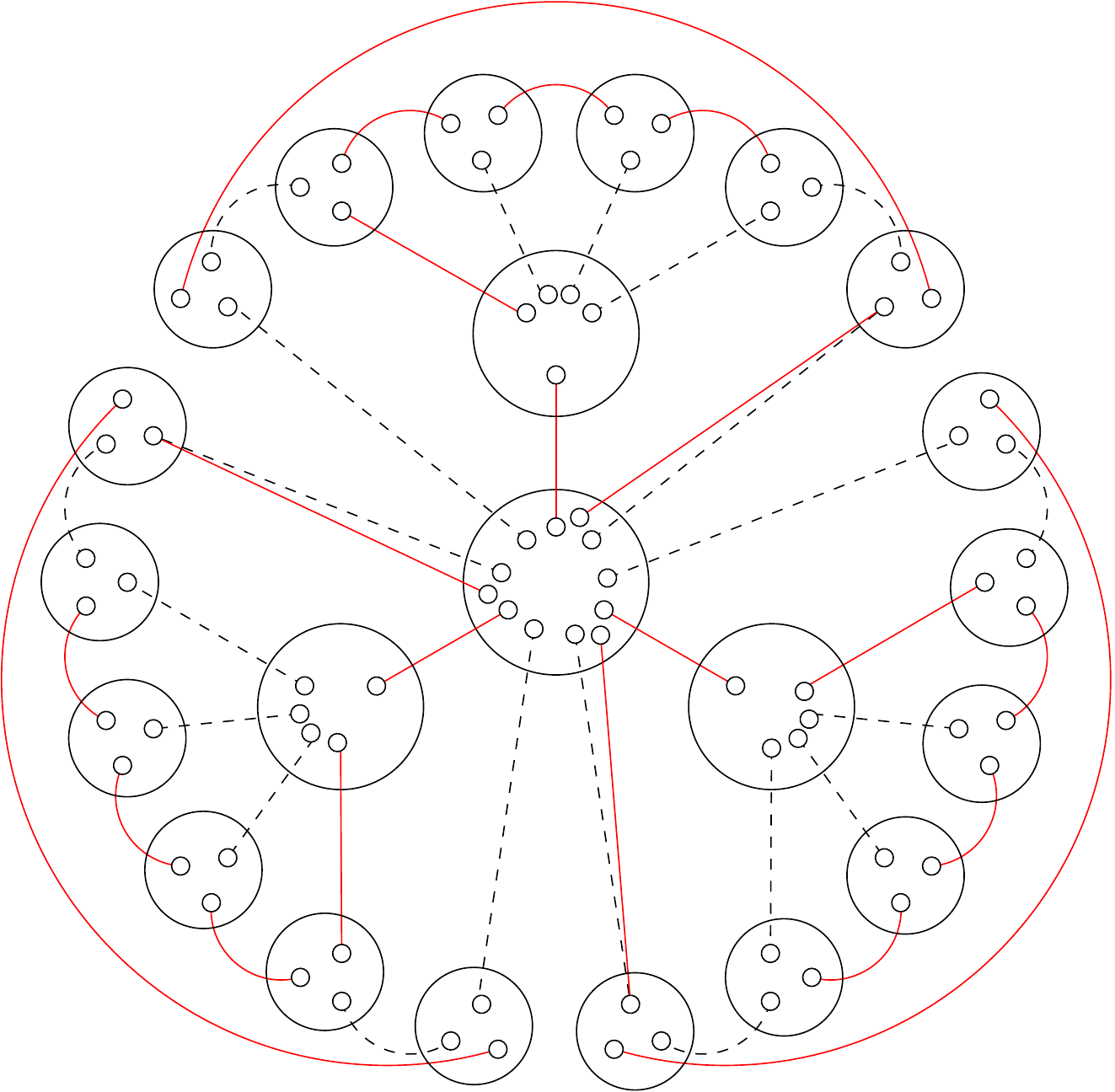}
\caption{\jbjr{Example of a $3$-connected 2T-graph $G$ such that the
    set of external edges almost form a matching and $G$ has no good
    ordering. The solid and dashed edges illustrate two spanning trees
    along the external edges which can be extended arbitrarily into
    the circuits. Note that there are 22 circuits and 42 external edges
    connecting these so all of these are needed by Theorem
    \ref{thm:tutte}. It also follows from Theorem \ref{thm:tutte}
    applied to the partition consisting of the seven circuits
    appearing from (roughly) 2 o'clock to 6 o'clock in the figure and
    the union of the remaining 15 circuits that the 4 external edges
    between these two collections are all needed and since they are
    incident with only 3 vertices of the seven circuits, there will be
    two external edges incident with the same vertex.  One gets
    further examples by enlarging the three paths on the rim of the
    figure.}}
\label{F7}
\end{center}
\end{figure}

Even if the graph has no conflict, then it is possible that it has no
good orientation. \jbjr{Indeed, using the example in Figure \ref{F6} we can now
  construct a complex example in Figure \ref{F7} below of a 2T-graph
  whose vertex set partitions into vertex sets of disjoint circuits so
  that $G$ has no good ordering. Note that it is necessary for the
  conclusion that there must be a global root in each of the three
  locally identical pieces of the graph that all of the circuits at
  the rim have exactly three vertices that are incident with external
  edges.}

\section{\jbj{2T-graphs which are disjoint unions of circuits}}
\label{decomposeintoGCsec}

In this section we consider 2T-graphs whose generic components are
circuits.  When we speak of a good orientation $D_{\prec}$ of a
2T-graph $G$, we use $s$ to denote the root of the out-branching $B^+$
and $t$ to denote the root of the in-branching $B^-$, where $B^+,B^-$
certify that $D$ is good (so $s$ is the first and $t$ is the last
vertex in the ordering $\prec$).

A circuit $H$ of $G$ is called a {\bf leaf} if there are exactly two
external edges between $H$ and some other circuit, that is, $H$
corresponds to a vertex in $\tilde{G}$ incident with two parallel
edges, otherwise $H$ is called {\bf internal}.

\jbj{\begin{theorem}
    \label{matchingcase}
    Let $G=(V,E)$ be a 2T-graph whose generic components are circuits.
    If the external edges in $G$ form a matching, then $G$ has a good
    ordering.
  \end{theorem}
  \pf Let $G_1,G_2,\ldots{},G_k$ be the circuits of $G$. We prove the
  theorem by by induction on $k$. When $k=1$, $G$ is itself a circuit
  and the result follows from Theorem \ref{thm:GChasgoodor}. So assume
  $k \geq 2$. Suppose first that some circuit $G_i$ is a leaf. By
  relabelling the circuits we may assume that $i=k$ and that $G_k$ 
  is connected to $G_{k-1}$ by a matching of 2 edges $uv,zw$, where
  $u,z\in V(G_{k-1})$ and $v,w\in V(G_k)$. By induction $G-V(G_k)$ has
  a good ordering ${\prec'}$. By renaming if necessary we can assume
  $u {\prec'} z$. By Theorem \ref{thm:GChasgoodor}, $G_k$ has a good
  ordering ${\prec''}$ such that $v$ is the first vertex and $w$ the
  last vertex of ${\prec''}$. Now we obtain a good ordering by
  inserting all the vertices of ${\prec''}$ just after $u$ in
  ${\prec'}$. Note that this corresponds to taking the union of the
  branchings $B_s^+,B_t^-$ that correspond to ${\prec'}$ and the
  branchings $\hat{B}_v^+,\hat{B}_w^-$ that correspond to ${\prec''}$
  by letting $A(\tilde{B}_s^+)=A(B_s^+)\cup A(\hat{B}_v^+)\cup\{uv\}$
  and $A(\tilde{B}_t^-)=A(B^-_t)\cup A(\hat{B}^-_w)\cup\{wz\}$.\\
  
  Suppose now that every $G_i$, $i\in [k]$ is internal. As $G$ and
  hence its quotient $\tilde{G}$ is a 2T-graph, there is a circuit
  $G_j$ such that there are exactly 3 edges $u_1v_1,u_2v_2,u_3v_3$,
  with $v_i\in V(G_j)$ connecting $V(G_j)$ to $V-V(G_j)$. Again we may
  assume that $j=k$. We may also assume w.l.o.g. that for some pair of
  spanning trees $T_1,T_2$ of $\tilde{G}$, the edges $u_1v_1,u_2v_2$
  belong to $T_1$ and $u_3v_3$ belongs to $T_2$ (so the vertex in
  $\tilde{G}$ corresponding to $G_k$ is a leaf in $T_2$). Note that
  this means that $u_1,u_2$ belong to different circuits
  $G_a,G_b$. Now let $H$ be obtained from $G$ by deleting the vertices
  of $V(G_k)$ and adding the edge $u_1u_2$. Then $V(H)$ decomposes
  into a disjoint union of vertex sets of circuits and set of edges
  connecting these form a matching. By induction there is a good
  ordering $\prec$ of $H$. Let $B^+_{s,0},B^-_{t,0}$ be a pair of
  arc-disjoint branchings that certify that $\prec$ is a good ordering
  of $H$. We are going to show how to insert the vertices of $V(G_k)$
  so that we obtain a good ordering of $G$. By renaming
  $u_1,u_2,v_1,v_2$ and possibly considering the reverse ordering
  $\stackrel{\leftarrow}{\prec}$ if necessary we can assume that $u_1
  {\prec} u_2$ and that the arc $u_1u_2$ belongs to $B^-_t$. We now
  consider the three possible positions of $u_3$ in the ordering
  ${\prec}$ (see Figure~\ref{F2}).

  \begin{itemize}
  \item $u_3{\prec}u_1{\prec}u_2$. By Theorem \ref{thm:GChasgoodor},
    $G_k$ has a good ordering ${\prec}_1$ of $G_k$ such that $v_3$ is
    the initial vertex and $v_2$ is the terminal vertex of
    ${\prec}_1$. Let $B^+_{v_3,1},B^-_{v_2,1}$ be arc-disjoint
    branchings (on $V(G_k)$) certifying that ${\prec}_1$ is good. Then
    we obtain a good ordering of $G$ by inserting all the vertices of
    ${\prec}_1$ just after $u_1$ in $\prec$ and we obtain the desired
    branchings $B^+_s,B^-_t$ by letting $A(B^+_s)=A(B^+_{s,0})\cup
    A(B^+_{v_3,1})\cup\{u_3v_3\}$ and
    $A(B^-_t)=A(B^-_{t,0}-u_1u_2)\cup
    A(B^-_{v_2,1})\cup\{u_1v_1,v_2u_2\}$.

  \item $u_1{\prec}u_2{\prec}u_3$. By Theorem \ref{thm:GChasgoodor},
    $G_k$ has a good ordering ${\prec}_2$ of $G_k$ such that $v_2$ is
    the initial vertex and $v_3$ is the terminal vertex of
    ${\prec}_2$. Let $B^+_{v_2,2},B^-_{v_3,2}$ be arc-disjoint
    branchings (on $V(G_k)$) certifying that ${\prec}_2$ is good. Then
    we obtain a good ordering of $G$ by inserting all the vertices of
    ${\prec}_2$ just after $u_2$ in $\prec$ and we obtain the desired
    branchings $B^+_s,B^-_t$ by letting $A(B^+_s)=A(B^+_{s,0})\cup
    A(B^+_{v_2,2})\cup\{u_2v_2\}$ and
    $A(B^-_t)=A(B^-_{t,0}-u_1u_2)\cup
    A(B^-_{v_3,2})\cup\{u_1v_1,v_3u_3\}$.

  \item $u_1{\prec}u_3{\prec}u_2$. %By Theorem \ref{thm:GChasgoodor}, $G_k$ has a good ordering ${\prec}_3$ of $G_k$ such that $v_3$ is the initial vertex and $v_2$ is the terminal vertex of ${\prec}_3$.
    Consider again the good ordering ${\prec}_1$ above and the branchings $B^+_{v_3,1},B^-_{v_2,1}$. Then we obtain a good ordering of $G$ by inserting all the vertices of ${\prec}_1$ just after $u_3$ in $\prec$ and we obtain the desired branchings $B^+_s,B^-_t$  by letting $A(B^+_s)=A(B^+_{s,0})\cup A(B^+_{v_3,1})\cup\{u_3v_3\}$ and $A(B^-_t)=A(B^-_{t,0}-u_1u_2)\cup A(B^-_{v_2,1})\cup\{u_1v_1,v_2u_2\}$.
  \end{itemize}

  As we saw, in all the possible cases we obtain a good ordering of
  $G$ together with a pair of arc-disjoint branchings which certify
  that the ordering is good so the proof is complete.  \qed

 \jbjb{ Figure \ref{F7} shows an example of a 2T-graph $G$ whose vertex set partitions into vertex sets of generic circuits such that the set of edges between different circuits almost forms a matching and the graph $G$ has no good ordering.}

\begin{figure}[h!]
\begin{center}
\scalebox{0.9}{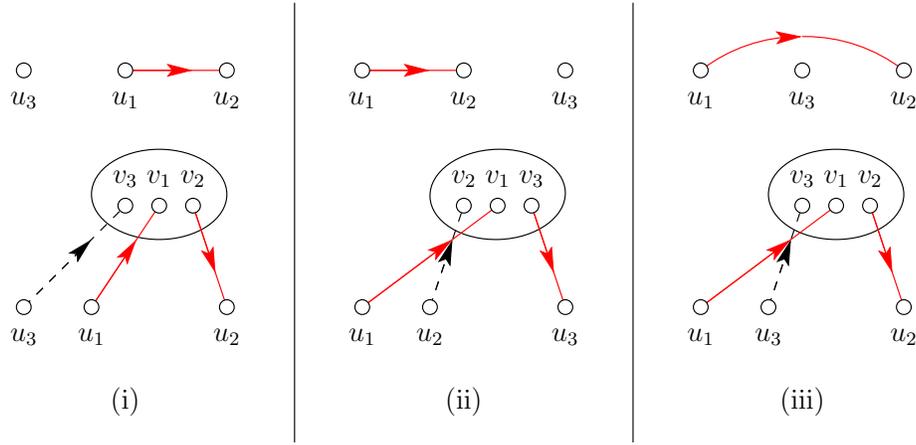}
\caption{How to lift a good ordering to a new circuit as in the proof of
  Theorem \ref{matchingcase}. In-branchings are displayed solid,
  out-branchings are dashed.  The first line displays the three
  possible orders of the relevant vertices (increasing left to right)
  as they occur in the proof. The second line displays the ordering of
  the augmented graph and how the branchings lead into and out of the
  new circuit; its local out- and in-root is the leftmost and
  rightmost $v_i$, respectively.}
\label{F2}
\end{center}
\end{figure}

}

A {\bf double tree} is any graph that one can obtain from a tree $T$
by adding one parallel edge for each edge of $T$. A {\bf double path}
is a double tree whose underlying simple graph is a path.

Recall that a subset $X\subset V$ with $2\leq |X|\leq |V|-2$ is
pendant at $x$ in $G$ if all edges between $X$ and $V(G)-X$ are
incident with $x$

\begin{definition}
  Let $G$ be a 2T-graph whose quotient graph is a double tree $T$.  An
  {\bf obstacle} in $G$ is a subgraph $H$ consisting of a subset of
  the circuits of $G$ and the edges between these such that the
  quotient graph of $G[V(H)]$ is a double path $T_H$ of $T$ so that
  \begin{itemize}
  \item $H$ contains circuits $C,C'$, possibly equal, and vertices
    $x\in C,y\in C'$, such that $x=y$ if $C=C'$ and there is an
    $(x,y)$-path $P$ in $H$ which uses only external edges of $H$ (so
    $P$ is also a path in $T_H$ between the two vertices corresponding
    to the circuits $C,C'$).
  \item $T-V(T_H)$ has at least two connected components $A,B$ and
    $V_A$ is pendant at $x$ and $V_B$ is pendant at $y$ in $G$, where
    $V_A$ (resp. $V_B$) is the union of those circuits of $G$ that
    correspond to the vertex set $A$ (resp. $B$) in $T$.
    \end{itemize}
\end{definition}

\begin{theorem}
\label{thm:doubleTchar}
Let $G$ be a 2T-graph whose quotient is a double tree $T$, then $G$
has a good ordering if and only if
\begin{itemize}
\item[(i)] $G$ has at most two pendant circuits and 
%\item[(ii)] $G$ has no internal circuit which is also pendant.
\item[(ii)] $G$ contains no obstacle.
\end{itemize}
\end{theorem}

\pf By Corollary \ref{cor:3pendant} we see that (i) must hold if $G$
has a good ordering.

Suppose that $G$ contains an obstacle $H$ but there exists a good
ordering ${\prec}$ with associated branchings $B^+_s,B^-_t$ in
$D=D_{\prec}$. Let $x,y$ be the special vertices according to the
definition. Suppose first that $x=y$ and let $C$ be the circuit that
contains $x$, let $v_C$ be the vertex of $T$ that corresponds to $C$
and let $V_A$, $V_B$ be the union of the vertex sets of circuits of
$G$ so that these correspond to distinct connected components $A,B$ of
$T-v_C$ and both $V_A$ and $V_B$ are pendant at $x$ in $G$.  By Lemma
\ref{lem:rootinpendant}, we may assume w.l.o.g that $s\in V_A$ and
$t\in V_B$. Then it is easy to see that $D$ contains two arcs
$a_1x,a_2x$ from $V_A$ to $x$ and the two arcs $xb_1,xb_2$ from $x$ to
$V_B$ and precisely one of the arcs $a_1x,a_2x$ is in $B^+_s$ and the
other is in $B^-_t$ and the same holds for the arcs $xb_1,xb_2$. Now
consider a vertex $z\in C-x$. The $(s,z)$-path in $B^+_s$ contains $x$
and the $(z,t)$-path in $B^-_t$ also contains $x$ so $D$ is not
acyclic, contradiction. Hence we must have $x\neq y$ and $x,y$ are in
different circuits (so $C\neq C'$). Again we let $V_A$ be the union of
vertices of circuits of $G$ so that $V_A$ is pendant at $x$ and
similarly let $V_B$ be the union of vertices of circuits of $G$ so
that $V_B$ is pendant at $y$. Again by Lemma \ref{lem:rootinpendant},
we may assume w.l.o.g that $s\in V_A$ and $t\in V_B$. As above $D$
must contain two arcs $a_1x,a_2x$ from $V_A$ to $x$ and the two arcs
$yb_1,yb_2$ from $y$ to $V_B$ and precisely one of the arcs
$a_1x,a_2x$ is in $B^+_s$ and the other is in $B^-_t$ and the same
holds for the arcs $yb_1,yb_2$. Let $C_1,C_2,\ldots{},C_r$, $r\geq 2$
be circuits of $G$ so that $C=C_1,C'=C_r$ and
$v_{C_1}v_{C_2}\ldots{}v_{C_r}$ is a path in $T$ which corresponds to
the $(x,y)$-path $P=x_1x_2\ldots{}x_r$, where $x=x_1,y=x_r$, that uses
only edges between different circuits in $G$ (by the definition of an
obstacle). As $s\in V_A$ and $t\in V_B$ the path $P$ must be a
directed $(x_1,x_r)$-path in $D$ and using that $D[V(C_1)]$ are
$D[V(C_r)]$ are acyclic we can conclude as above that the arc $x_1x_2$
is an arc of $B^+_s$ and the arc $x_{r-1}x_r$ is an arc of
$B^-_t$ \jbjb{(if $x_1x_2$ would not be  an arc of $B^+_s$, then $D[V(C_1)]$ would contain a directed path from $x_1$ to the end vertex $z$ of the other arc leaving $V(C_1)$ and also a directed path from $z$ to $x_1$, implying that $D[V(C_1)]$ would not be acyclic)}. Thus it follows that for some index $1<j<r$ the arc
$x_{j-1}x_j$ is an arc of $B^+_s$ and the arc $x_jx_{j+1}$ is an arc
of $B^-_t$. This implies that for every $z\in C_j$ the $(s,z)$-path of
$B^+_s$ and the $(z,t)$-path of $B^-_t$ contains $x_j$, contradicting
that $D$ is acyclic.

Suppose now that $G$ satisfies (i) and (ii). We shall prove by
induction on the number, $k$, of circuits in $G$ that $G$ has a good
orientation. The base case $k=1$ follows from Theorem
\ref{thm:GChasgoodor} so we may proceed to the induction step.\\

Suppose first that $G$ has a leaf circuit $G_h$ that is not
pendant. Let $v_{h'}$ be the neighbour of $v_h$ in $\tilde{G}$ and let
$G_{h'}$ be the circuit of $G$ corresponding to $v_{h'}$. As $G_h$ is
not pendant the two edges $xx',zz'$ between $G_h$ and $G_{h'}$ have
distinct end vertices $x,z$ in $V(G_h)$ and distinct end vertices
$x',z'$ in $V(G_{h'})$.  By induction $G-G_h$ has a good orientation
$D'$ and we may assume, by reversing all arcs, if necessary, that $x'$
occurs before $z'$ in the ordering ${\prec'}$ that induces $D'$. By
Theorem \ref{thm:GChasgoodor}, $G_h$ has a good orientation $D''$
where $x$ is the out-root and $z$ is the in-root. Now we obtain a good
orientation $D$ by adding the two arcs $x'x$ and $zz'$ and
using the first in the out-branching rooted at $x$ and the later in
the in-branching. \\

Thus we can assume from now on that every leaf component of $G$ is
pendant and now it follows from Corollary \ref{cor:3pendant} that $G$
is a double path whose circuits we can assume are ordered as
$G_1,G_2,\ldots{},G_k$ in the ordering that the corresponding vertices
$v_1,v_2,\ldots{},v_k$ appear in the quotient $\tilde{G}$.\\

We prove the following stronger statement which will imply that $G$
has a good orientation.

\begin{claim}
  \label{claim1}
  Let $G$ be a double path having no obstacle and whose circuits are
  ordered as $G_1,G_2,\ldots{},G_k$. Let $s\in V(G_1)$ be any
  vertex, except $a$ if $G_1$ is pendant at $a\in V(G_1)$ and let
  $t\in V(G_k)$ be any vertex except $b$ if $G_k$ is pendant at $b\in
  V(G_k)$ (such vertices are called {\bf candidates for roots}). Then
  has a good orientation $D_{\prec}$ so that $s$ is the first vertex
  (root of the out-branching) and $t$ is the last vertex in $\prec$ if
  and only if none of the following hold.
\begin{enumerate}

\item[(a)] There is an $(s,t)$-path $P$ in $G$ which uses only
  external edges. %no edge which is inside some $G_i$.
\item[(b)] $t$ is an end vertex of one of the edges from $G_{k-1}$,
  there is an index $i\in [k-1]$ so that the two edges from $G_i$ to
  $G_{i+1}$ are incident with the same vertex $x$ of $G_{i+1}$ and
  there is an $(x,t)$-path in $G$ which uses only external edges.
  %no edge which is inside some $G_i$.
\item[(c)] $s$ is an end vertex of one of the edges from $G_1$ to
  $G_2$, there is an index $j\in [k]\setminus \{1\}$ so that the two edges from
  $G_{j-1}$ to $G_j$ are incident with the same vertex $y$ of
  $G_{j-1}$ and there is an $(s,y)$-path in $G$ which uses only external
  edges.
% no edge which is inside some $G_q$.
  %\item[(d)] There exist indices $1<i<j<$ and vertices $u\in V(G_i),v\in V(G_j)$ so that $u$ is incident with both edges from $V(G_{i-1})$ to $V(G_i)$, $v$ is incident with both edges from $V(G_j)$ to $V(G_{j+1})$ and $G$ contains a $(u,v)$-path which uses no edge which is inside some $G_r$.
\end{enumerate}
\end{claim}

\noindent{}{\bf Proof of claim:} Note that if
${\prec}:\ v_1,\ldots{},v_n$ is a good ordering with $s=v_1$ and
$t=v_n$, then, in the corresponding acyclic orientation $D_{\prec}$
the two edges between $G_i$ and $G_{i+1}$ are both oriented towards
$G_{i+1}$ and for every pair of arc-disjoint branchings $B^+_s,B^-_t$
in $D$, exactly one of these arcs belong to $B^+_s$ and the other to
$B^-_t$.\\

We first show that if $G,s,t$ satisfy any of (a)-(c), then there is no
good ordering $v_1,\ldots{},v_n$ with $s=v_1$ and $t=v_n$.\\

Suppose that $G$ has a good ordering $v_1,\ldots{},v_n$ with $s=v_1$
and $t=v_n$ and let $B^+_s,B^-_t$ be a pair of arc-disjoint branchings
in the acyclic digraph $D=D_{\prec}$.

If (a) holds, then let $P=x_1x_2\ldots{}x_k$ be a path from $s=x_1$ to
$t=x_k$ so that each edge $x_ix_{i+1}$, $i\in [k-1]$ has one end
vertex in $G_i$ and the other in $G_{i+1}$. As $B^+_s$ induces and
out-branching from $s$ in the acyclic digraph $D[V(G_1)]$, we must
have that the arc $sx_2$ belongs to $B^+_s$.
% This means that the other arc between $V(G_1)$ and $V(G_2)$ belongs to $B^__t$
By a similar argument, the arc $x_{k-1}t$ belongs to $B^-_t$. Hence
there is an index $1<i<k$ such that the arc $x_{i-1}x_i$ is in $B^+_s$
and the arc $x_ix_{i+1}$ is in $B^-_t$. However this implies that
$x_i$ is both an out-root and an in-root in $D[V(G_i)]$, contradicting
that $D$ is acyclic. So (a) cannot hold if there is a good
ordering.\\ If (b) holds, then let \jbjb{$x\in V(G_{i+1})$} be the vertex
incident with both edges between \jbjb{$G_{i}$ and $G_{i+1}$} and let
\jbjb{$x=x_{i+1}x_{i+2}\ldots{},x_{k-1}t$} be an $(x,t)$-path in $G$ so that
each edge $x_jx_{j+1}$, $\jbjb{i+1}\leq j\leq k$ has one end vertex in \jbjb{$G_j$}
and the other in \jbjb{$G_{j+1}$}. As above, we conclude that the arc
$x_{k-1}t$ belongs to $B^-_s$ and that there exists an index $j$ with
$i\leq j$ so that $x_j$ is the head an arc of $B^+_s$ coming from
$G_{j-1}$ and the tail of an arc of $B^-_t$ going to $G_{j+1}$. As
above this again contradicts that $D$ is acyclic with branchings
$B^+_s,B^-_t$.  Analogously we see that (c) cannot hold when there is
a good ordering.  \iffalse Finally suppose that (d) holds and let
$u,v$ be vertices and let $Q$ be the path between them as in
(d). Since $u$ is the head of an arc from each of the branchings and
$v$ is the tail of an arc from both it is easy to conclude that there
is a vertex $w\in G_j$ on $Q$ which is the head of an arc from
$G_{j-1}$ that belongs to $B^+_s$ and also tail of an arc from $G_j$
to $G_{j+1}$ which belongs to $B^-_t$. As before this contradicts that
$D$ is acyclic with branchings $B^+_s,B^-_t$.\\ \fi

Suppose now that none of (a)-(c) hold.  We prove the existence of a
good orientation by induction on $k$. For $k=1$ the claim follows from
Theorem \ref{thm:GChasgoodor}.

Suppose next that $k=2$.
\jbjb{Let $u_1u_2$ and $v_1,v_2$ with $u_1,v_1\in V(G_1)$ be
the two edges between $G_1$ and $G_2$. Suppose first that
$t\not\in \{u_2,v_2\}$. Since $G_1$ is not pendant at
$s$, we can assume w.l.o.g. that $s\neq v_1$.  By Theorem \ref{thm:GChasgoodor}, there is a
good orientation of $G_1$ in which $s$ is the out-root and $v_1$ is the
in-root and a good orientation of $G_2$ in which $u_2$ is the out-root
and $t$ is the in-root. Thus we obtain the desired orientation by
adding the arc $u_1u_2$ to the union of the two out-branchings and
the arc $v_1v_2$ to the union of the two in-branchings. Suppose now that $t\in \{u_2,v_2\}$. Without loss of generality $t=u_2$. Since (a) does not hold, we know that $s\neq u_1$. By
Theorem \ref{thm:GChasgoodor}, $G_1$ has a good orientation with $s$
and out-root and $u_1$ as in-root and $G_2$ has a good orientation with
$v_2$ as out-root and $t$ as in-root. Now we obtain the desired
branchings by adding the arc $u_1u_2$ to the union of the two in-branchings
and the arc $v_1v_2$ to the union of the two out-branchings.
}

Assume that $k\geq 3$ and that the claim holds for all double paths which satisfy none of (a)-(c) and have  fewer than $k$ circuits. Let $s\in V(G_1),t\in V(G_k)$ be candidates for roots and let $xx',zz'$ be the two edges between $G_1$ and $G_2$. \jbjb{Without loss of generality we have $s\neq z$.}
% As $G$ has no internal pendant circuit at least one of $x',z'$, say w.l.o.g. $x'$, does not cover all edges between $G_2$ and $G_3$.
Note that (b) cannot  hold for $s',t$ in $G'=G-G_1$ when $s'\in\{x',z'\}$ because $G'$ is an induced subgraph of $G$. Suppose that (a) holds for $(G',x',t)$. Then $z'\neq x'$ as (b) does not hold for $G$. Now (a) cannot hold for $(G',z',t)$ as this would imply that (b) holds in $G$. For the same reason (c) cannot  hold for $(G',z',t)$. Thus if (a) holds for $(G',x',t)$, then none of (a)-(c) hold for $(G',z',t)$. If (c) holds for $(G',x',t)$ we conclude that none of (a),(c) hold for $(G',z',t)$, because both would imply that $G$ contains an obstacle. Let $s'=x'$ unless one of (a)-(c) holds for $x'$ and in that case \jbjb{$s\neq x$ must hold and } we let $s'=z'$. By the arguments above, none of (a)-(c) hold for $(G',s',t)$.\\

By induction $G'$ has a good orientation where $s'$ is the out-root
and $t$ is the in-root and by Theorem \ref{thm:GChasgoodor}, $G_1$ has
a good orientation in which $s$ is the out-root and $z$ is the
in-root. Let $a$ be the arc $xx'$ if $s'=x'$ and otherwise let $a$ be
the arc $zz'$. Now adding $a$ to the union of the two out-branchings
and the other arc from $G_1$ to $G_2$ to the union of the two
in-branchings, we obtain the desired good orientation. This completes
the proof of Claim \ref{claim1}.\\

Now we are ready to conclude the proof of Theorem
\ref{thm:doubleTchar}. As $G_1,G_k$ are circuits they both have at
least 2 vertices. If we can choose $s\in V(G_1)$ and $t\in V(G_k)$ so
that none of these two vertices are incident with edges to the other
circuits, then we are done by the Claim \ref{claim1}, so either
$|V(G_1)|=2$ or $|V(G_k)|=2$ or both. Suppose w.l.o.g. that
$|V(G_1)|=2$ and that the two edges from $G_1$ to $V-G_1$ are incident
with different vertices $u,v\in V(G_1)$. As $V(G_1)$ is pendant, these
two edges end in the same vertex $x$. If we can choose $t\in V(G_k)$
so that it is not incident with any of the edges between $G_{k-1}$ and
$G_k$, then we are done, so we may assume that we also have
$V(G_k)=\{z,w\}$ and that the edges between $G_{k-1}$ and $G_k$ are
$yz,yw$ for some $y\in V(G_{k-1})$ ($G_k$ is pendant).  Now it follows
from the fact that (ii) holds that every $(x,y)$-path in $G$ uses an
edge which is inside some $G_i$ we can take $s$ and $t$ freely among
$u,v$, respectively $z,w$ and conclude by the claim (none of (a)-(c)
can hold).  \qed

\iffalse
we just need to recall that each circuit has at least two vertices and it is easy to check that  we can choose $s\in V(G_1), t\in V(G_k)$ so that none of (a)-(c) hold.
\fi

\iffalse
 Let $x\in V(G_1)$ be the vertex incident with the two edges $xy,xz$ from $G_1$ to $G_2$. By induction $G-G_1$ has a good ordering ${\prec}^1$  in which $y$ is the out-root of $D_{{\prec}^1})$ and some vertex $t\in V(G_k)$ is the in-root (if $k=2$ we can choose $t\neq y$). By Theorem \ref{thm:GChasgoodor}, $G_1$ has a good ordering 
${\prec}^2$ such that $s$ is the out-root and $x$ is the in-root in $D_{{\prec}^2})$, where $s$ is an arbitrary vertex of $V(G_1)-x$. Thus we obtain the desired good orientation of $D$ by adding the arcs $x\dom y$ and $x\dom z$ and using the first in the out-branching rooted at $s$ and the later in the in-branching rooted at $t$. 
\qed

 \fi

\section{Remarks and open problems}\label{remarksec}

Let us start by recalling that the following is an immediate
consequence of Theorem \ref{thm:GChasgoodor} as we first find a good
ordering of the circuit and then orient the remaining edges according to
that ordering.

\begin{corollary}
  \label{cor:containsGC}
Every graph which contains  a circuit as a spanning subgraph
has a good ordering.
\end{corollary}

\begin{conjecture}
There exists a polynomial algorithm for deciding whether a 2T-graph
has a good ordering.
\end{conjecture}

\begin{problem}
What is the complexity of deciding whether a given graph has a good ordering?
\end{problem}

Two of the authors of the current paper proved the following generalization of Theorem \ref{4reg4con}. Note that its proof is more complicated than that of Theorem  \ref{4reg4con}.

\begin{theorem}\cite{bang4reg4con}
\label{all4reg4congood}
Every 4-regular 4-connected graph has a good orientation.
\end{theorem}

\iffalse
\begin{problem}
What is the complexity of deciding whether a mixed graph $M=(V,E\cup
A)$ which is acyclic (no directed cycle in $D=(V,A)$) and whose
underlying graph is a circuit, has a good ordering that is consistent
with the arcs in $A$?
\end{problem}

\begin{problem}
What is the complexity of deciding whether a mixed graph $M=(V,E\cup
A)$ which is acyclic (no directed cycle in $D=(V,A)$) and whose
underlying graph is a 2T-graph, has a good ordering that is consistent
with the arcs in $A$?
\end{problem}
\fi

  Let $D=(V,A)$ be a digraph and let $s,t$ be distinct vertices of
  $V$. An {\bf $(s,t)$-ordering} of $D$ is an ordering
  ${\prec}:\ v_1,v_2,\ldots{},v_n$ with $v_1=s,v_n=t$ such that every
  vertex $v_i$ with $i<n$ has an arc to some $v_j$ with $i<j$ and
  every vertex $v_r$ with $r>1$ has an arc from some $v_p$ with
  $p<r$. It is easy to see that $D$ has such an ordering if and only
  it it has a spanning acyclic digraph in with branchings
  $B^+_s,B^-_t$. These branchings are not necessarily arc-disjoint but
  it is clear that if $D$ has a good ordering with $s$ as the initial
  and $t$ the terminal vertex then this ordering is also an
  $(s,t)$-ordering. Hence having an $(s,t)$-ordering is a necessary
  condition for having a good ordering with $s$ as the initial and
  $t$ as the terminal vertex.

\begin{theorem}
  \label{thm:storderNPC}
  It is NP-complete to decide whether a digraph $D=(V,A)$ with
  prescribed vertices $s,t\in V$ has an $(s,t)$-ordering.
\end{theorem}
\pf The so called {\sc betweenness} problem is as follows: given a set
$S$ and a collection of triples $(x_i,y_i,z_i)$, $i\in [m]$,
consisting of three distinct elements of $S$; is there a total order
on $S$ (called a betweenness order on $S$) so that for each of the
triples we have either $x_i<y_i<z_i$ or $z_i<y_i<x_i$? {\sc
  Betweenness} is NP-complete \cite{opatrnySJC8}. Given an instance
$[S, (x_i,y_i,z_i), i\in [m]]$ of {\sc betweenness} we construct the
following digraph $D$. The vertex set $V$ of $D$ is constructed as
follows: first take $5m$ vertices
$$a_1,\ldots{},a_m,b_1,\ldots{},b_m,c_1,\ldots{},c_m,d_1,\ldots{},d_m,e_1,\ldots{},e_m$$
where $\{a_i,b_i,c_i\}$ corresponds to the triple $(x_i,y_i,z_i)$ and
then identify those vertices in the set
$\{a_1,\ldots{},a_m,b_1,\ldots{},b_m,c_1,\ldots{},c_m\}$ that
correspond to the same element of $S$. Then, add two more vertices:
$s$ and $t$. The arc set of $D$ consists of an arc from $s$ to each
vertex of $\{a_1,\ldots{},a_m,c_1,\ldots{},c_m\}$, an arc from each
such vertex to $t$ and the following $6m$ arcs which model the
betweenness conditions: for each triple $(x_i,y_i,z_i)$ $D$ contains
the arcs $a_id_i,c_id_i,d_ib_i,$ $b_ie_i,e_ia_i,e_ic_i$. \jbjb{Clearly $D$ can be constructed in polynomial time.} We claim that
$D$ has an $(s,t)$-ordering if and only if there is a betweenness
total ordering of $S$. Suppose first that $D$ has an
$(s,t)$-ordering. The vertices $d_i,b_i,e_i$ must occur in that order
as $b_i$ is the unique out-neighbour (in-neighbour) of $d_i$
($e_i$). As each $a_i,c_i$ are the only in-neighbours (out-neighbours)
of $d_i$ ($e_i$) in $D$ the vertices $a_i,c_i$ cannot both occur after
(before) $d_i$ ($e_i$) so the vertices in $\{a_i,b_i,c_i\}$ will occur
either in the order $a_i,b_i,c_i$ or in the order $c_i,b_i,a_i$. Thus
taking the same order for the elements in $S$ as for the corresponding
vertices of $D$, we obtain a betweenness total order. Conversely, if
we are given a betweenness total order for $S$ we just place the
vertices in $\{a_1,\ldots{},a_m,b_1,\ldots{},b_m,c_1,\ldots{},c_m\}$
in the order that the corresponding elements of $S$ occur and then
insert each vertex $d_i$ ($e_i$) anywhere between $a_i$ and $b_i$
($b_i$ and $c_i$) if the triple $(x_i,y_i,z_i)$ is ordered as
$x_i<y_i<z_i$ and otherwise we insert $d_i$ ($e_i$) anywhere between
$c_i$ and $b_i$ ($b_i$ and $a_i$). Finally insert $s$ as the first
element and $t$ as the last element.  Now every vertex different from
$s,t$ has an earlier in-neighbour and a later out-neighbour, so it is
an $(s,t)$-ordering. \qed

If $D$ is semicomplete digraph, that is, a digraph with no pair of non-adjacent vertices, then $D$ has an $(s,t)$-ordering for a given pair of distinct vertices $s,t$ if and only if $D$ has a Hamiltonian path from $s$ to $t$ \cite{bangJGT20a,thomassenJCT28}. It was shown in \cite{bangJA13} that there exists a polynomial algorithm for deciding the existence of such a path in a given semicomplete digraph so for semicomplete digraphs the $(s,t)$-ordering problem is polynomially solvable.

\iffalse
\begin{question}[For us!]
  What is the complexity of deciding whether a strong digraph $D=(V,A)$ has an $(s,t)$-ordering for some choice of distinct vertices $s,t\in V$?
\end{question}\fi

\begin{corollary}
  It is NP-complete to decide if a strong digraph $D=(V,A)$ has an $(p,q)$-ordering for some choice of distinct vertices $p,q\in V$
\end{corollary}

\pf Let $D'$ be the digraph that we obtain from the digraph $D$ in the proof above by adding the arc $ts$. Then $D'$ is strong and it is easy to see that the only possible pair for which there could exists a $(p,q)$-ordering is the pair $p=s,q=t$: for each triple $(x_i,y_i,z_i)$ the corresponding vertices in $D$ must occur either in the order $a_i,d_i,b_i,e_i,c_i$ or in the order
$c_i,d_i,b_i,e_i,a_i$ and in both cases $s$ must be before all these vertices and $t$ must be after all these vertices.\qed

 \begin{problem}
  What is the complexity of deciding whether a digraph which has a pair of arc-disjoint branchings $B^+_s,B^-_t$ has such a pair whose union (of the arcs) is an acyclic digraph?
\end{problem}

\bibliography{refs}

\begin{thebibliography}{10}

\bibitem{bangJCT51}
J.~Bang-Jensen.
\newblock {Edge-disjoint in- and out-branchings in tournaments and related path
  problems}.
\newblock {\em J. Combin. Theory Ser. B}, 51(1):1--23, 1991.

\bibitem{bangJGT20a}
J.~Bang-Jensen.
\newblock {Digraphs with the path-merging property}.
\newblock {\em J. Graph Theory}, 20(2):255--265, 1995.

\bibitem{bang2009}
J.~Bang-Jensen and G.~Gutin.
\newblock {\em {Digraphs: Theory, Algorithms and Applications}}.
\newblock Springer-Verlag, London, 2nd edition, 2009.

\bibitem{bangJCT102}
J.~Bang-Jensen and J.~Huang.
\newblock {Decomposing locally semicomplete digraphs into strong spanning
  subdigraphs}.
\newblock {\em J. Combin. Theory Ser. B}, 102:701--714, 2010.

\bibitem{bang4reg4con}
J~Bang-Jensen and M.~Kriesell.
\newblock Good acyclic orientations of 4-regular 4-connected graphs.
\newblock {\em manuscript}, 2019.

\bibitem{bangJA13}
J.~Bang-Jensen, Y.~Manoussakis, and C.~Thomassen.
\newblock {A polynomial algorithm for Hamiltonian-connectedness in semicomplete
  digraphs.}
\newblock {\em J. Algor.}, 13(1):114--127, 1992.

\bibitem{bangDAM161a}
J.~Bang-Jensen and S.~Simonsen.
\newblock {Arc-disjoint paths and trees in 2-regular digraphs}.
\newblock {\em Discrete Appl. Math.}, 161:2724--2730, 2013.

\bibitem{bangJGT42}
J.~Bang-Jensen, S.~Thomass{\'e}, and A.~Yeo.
\newblock {Small degree out-branchings}.
\newblock {\em J. Graph Theory}, 42(4):297--307, 2003.

\bibitem{bangC24}
J.~Bang-Jensen and A.~Yeo.
\newblock {Decomposing {$k$}-arc-strong tournaments into strong spanning
  subdigraphs}.
\newblock {\em Combinatorica}, 24(3):331--349, 2004.

\bibitem{bangTCS438}
J.~Bang-Jensen and A.~Yeo.
\newblock {Arc-disjoint spanning sub(di)graphs in digraphs}.
\newblock {\em Theor. Comput. Sci.}, 438:48--54, 2012.

\bibitem{bergJCT88}
A.~R. Berg and T.~Jord{\'{a}}n.
\newblock {A proof of Connelly's conjecture on 3-connected circuits of the
  rigidity matroid}.
\newblock {\em J. Comb. Theory, Ser. {B}}, 88(1):77--97, 2003.

\bibitem{bergLNCS2832}
A.R. Berg and T.~Jord{\'{a}}n.
\newblock {Algorithms for Graph Rigidity and Scene Analysis}.
\newblock In {\em Algorithms - {ESA} 2003, 11th Annual European Symposium,
  Budapest, Hungary, September 16-19, 2003, Proceedings}, pages 78--89, 2003.

\bibitem{edmonds1973}
J.~Edmonds.
\newblock {Edge-disjoint branchings}.
\newblock In {\em {Combinatorial Algorithms}}, pages 91--96. Academic Press,
  1973.

\bibitem{garey1979}
M.R. Garey and D.S. Johnson.
\newblock {\em {Computers and intractability}}.
\newblock W. H. Freeman, San Francisco, 1979.

\bibitem{kocholJCT78}
M.~Kochol.
\newblock {Equivalence of Fleischner's and Thomassen's Conjectures}.
\newblock {\em J. Comb. Theory, Ser. {B}}, 78(2):277--279, 2000.

\bibitem{lamanJEM4}
G.~Laman.
\newblock On graphs and rigidity of plane skeletal structures.
\newblock {\em J. Eng. Math.}, 4:331--340, 1970.

\bibitem{matthewsJGT8}
M.~M. Matthews and D.~P. Sumner.
\newblock {Hamiltonian results in $K_{1,3}$-free graphs}.
\newblock {\em Journal of Graph Theory}, 8(1):139--146, 1984.

\bibitem{nashwilliamsJLMS39}
C.St.J.A. Nash-Williams.
\newblock {Decomposition of finite graphs into forests}.
\newblock {\em J. London Math. Soc.}, 39:12, 1964.

\bibitem{opatrnySJC8}
J.~Opatrny.
\newblock Total ordering problem.
\newblock {\em Siam J. Comput.}, 8:111--114, 1979.

\bibitem{recski1989}
A.~Recski.
\newblock {\em {Matroid theory and its applications in electric network theory
  and in statics}}.
\newblock Springer-Verlag, Berlin, 1989.

\bibitem{thomassenJCT28}
C.~Thomassen.
\newblock {Hamiltonian-connected tournaments}.
\newblock {\em J. Combin. Theory Ser. B}, 28(2):142--163, 1980.

\bibitem{thomassenJGT10}
C.~Thomassen.
\newblock {Reflections on graph theory}.
\newblock {\em J. Graph Theory}, 10(3):309--324, 1986.

\bibitem{tutteJLMS36}
W.T. Tutte.
\newblock {On the problem of decomposing a graph into {$n$} connected factors}.
\newblock {\em J. London Math. Soc.}, 36:221--230, 1961.

\bibitem{ryjacekJCT70}
Ryj{\'{a}}cek Z.
\newblock {On a Closure Concept in Claw-Free Graphs}.
\newblock {\em J. Comb. Theory, Ser. {B}}, 70(2):217--224, 1997.

\end{thebibliography}

\end{document}